\documentclass[aos,MSNbibl,nameyear,dvips]{arximspdf}
\usepackage{mathbh}
\usepackage{graphicx}

\doi{10.1214/14-AOS1240} 
\volume{42}
\issue{5}
\pubyear{2014}
\firstpage{1879}
\lastpage{1910}
\docsubty{FLA}

\makeatletter
\newcommand{\MSE}{\operatorname{MSE}}
\def\mathds{\mathbh}
\newcommand{\eqref}[1]{(\ref{#1})}
\newproclaim{defn}{Definition}[section]
\newtheorem{lem}{Lemma}[section]
\newtheorem{thmm}{Theorem}[section]
\newtheorem{prop}{Proposition}[section]
\newtheorem{cor}{Corollary}[section]

\newproclaim{rk}{Remark}

\newcommand{\1}{\mathds{1}}
\newcommand{\st}{| }
\newcommand{\iid}{i.i.d.}
\newcommand{\argmin}{\mathop{\operatorname{Argmin}}}

\newcommand{\R}{\mathbb{R}}
\newcommand{\N}{\mathbb{N}}
\newcommand{\E}{\mathbb{E}}
\renewcommand{\P}{\mathbb{P}}
\newcommand{\Var}{\operatorname{Var}}
\newcommand{\Cov}{\operatorname{Cov}}
\newcommand{\B}{\mathbb{B}}
\newcommand{\crit}{\operatorname{crit}}
\newcommand{\card}{\operatorname{Card}} 
\newcommand{\sh}{\widehat{s}}
\newcommand{\mh}{\widehat{m}}
\newcommand{\M}{\mathcal{M}}
\newcommand{\Z}{\mathcal{Z}}
\renewcommand{\S}{\mathcal{S}}
\newcommand{\Pg}{P\gamma}
\newcommand{\g}{\gamma}
\newcommand{\Rh}{\widehat{R}}
\newcommand{\bayes}{s}

\makeatother

\begin{document}
\begin{frontmatter}

\title{Optimal cross-validation in density estimation with~the $L^2$-loss\thanksref{T2}}
\runtitle{Optimal cross-validation}
\thankstext{T2}{Supported by the French
Agence Nationale de la Recherche (ANR) under reference
ANR-09-JCJC-0027-01
and ANR-11-BS01-0010.}

\begin{aug}
\author{\fnms{Alain} \snm{Celisse}\corref{}\ead[label=e1]{celisse@math.univ-lille1.fr}}
\runauthor{A. Celisse}
\address{Laboratoire de Math\'ematiques Painlev\'e\\
UMR 8524 CNRS-Universit\'e Lille 1\\
\textsc{Modal} project-team INRIA\\
F-59655, Villeneuve d'Ascq Cedex\\
France\\
\printead{e1}}
\affiliation{UMR 8524 CNRS--Universit\'e Lille 1}
\end{aug}

\received{\smonth{3} \syear{2012}}
\revised{\smonth{5} \syear{2014}}

%
\begin{abstract}
We analyze the performance of cross-validation (CV) in the density
estimation framework with two purposes: (i) \emph{risk estimation}
and (ii) \emph{model selection}.
The main focus is given to the so-called leave-$p$-out CV procedure
(Lpo), where $p$ denotes the cardinality of the test set.
Closed-form expressions are settled for the Lpo estimator of the risk
of projection estimators. These expressions provide a great improvement
upon $V$-fold cross-validation in terms of variability and computational
complexity.

From a theoretical point of view, closed-form expressions also enable
to study the Lpo performance in terms of risk estimation.
The optimality of leave-one-out (Loo), that is Lpo with $p=1$, is
proved among CV procedures used for risk estimation.
Two model selection frameworks are also considered: \emph{estimation},
as opposed to \emph{identification}.
For estimation with finite sample size $n$, optimality is achieved for
$p$ large enough [with $p/n =o(1)$] to balance the overfitting
resulting from the structure of the model collection.
For identification, model selection consistency is settled for Lpo as
long as $p/n$ is conveniently related to the rate of convergence of the
best estimator in the collection: (i) $p/n\to1$ as $n\to+\infty$
with a parametric rate, and (ii)~$p/n=o(1)$ with some nonparametric
estimators.
These theoretical results are validated by simulation experiments.
\end{abstract}

%

%
\end{frontmatter}

\section{Introduction}


%
For estimating a target quantity denoted by $s$, let\break  $\{ S_m \}_{m\in
\M
}$ denote a collection of sets of candidate parameters indexed by $\M$.
From each $S_m$ called a \emph{model}, an estimator $\widehat{s}_m$
of $s$ is computed.
The goal of model selection is to design a criterion $\crit\dvtx \M
\to\R
^+$ such that minimizing $\crit(\cdot)$ over $\M$ provides a final
estimator $\sh_{\mh}$ that is ``optimal.''
Among various strategies of model selection, \emph{model selection via
penalization} has been introduced in the seminal
papers by \citet{Mall73,Akai73,Schw78} on, respectively, AIC, $C_p$ and
BIC criteria.
However, since AIC and BIC are derived from asymptotic arguments, their
performances crucially depend on model collection and sample size
[see \citet{BaGH09}].\vadjust{\goodbreak}

More recently, \citeauthor{BiMa97} (\citeyear{BiMa97,BiMa01,BiMa06}) have
developed a
nonasymptotic approach inspired from the pioneering work of
\citet{BaCo91}.
It relies on concentration inequalities [\citet{Tala96,Led:2001}] and
aims at deriving \emph{oracle inequalities} such as
%
\begin{equation}
\label{exp.oracle.inequality.general} \ell( s,\sh_{\mh} ) \leq C \inf_{m\in\M}
\bigl\{ \ell( s,\widehat {s}_m ) \bigr\} + r_n
\end{equation}
with probability larger than $1-c/n^2$, where $c>0$ is a constant,
$\ell
(s,t)$ is a measure of the gap between parameters $s$ and $t$, $r_n$ is
a remainder term with respect to $\inf_m \ell( s, \widehat{s}_m )$, and
$C\geq1$ denotes a constant independent of $s$.
The closer $C$ to~1 and the smaller $r_n$, the better the model
selection procedure.
If $C=C_n\to1$ as $n\to+\infty$, the model selection procedure is
said \emph{asymptotically optimal} (or efficient)
[see, e.g., \citet{ArCe2010survey}].
Note that other asymptotic optimality properties have been studied in
the literature. For instance, a model selection procedure satisfying
\[
\P[ \mh= m_0 ] \mathop{\longrightarrow}_{n\to+\infty} 1,
\]
where $m_0$ denotes a fixed given model is said \emph{model selection
consistent} [see \citet{Sha:1997} for a study of various model selection
procedures in terms of model selection consistency].

In the density estimation framework, model selection with deterministic
penalties has been addressed: (i) for Kullback--Leibler divergence by
\citet{BarBirMas:1999},
\citeauthor{Cast99} (\citeyear{Cast99,Cast03}),
\citet{YangBarron1998} and further
studied in
\citet{BiRo06}, and (ii) for quadratic risk and projection estimators by
\citet{BiMa97} and \citet{BarBirMas:1999}.

The aforementioned approaches rely on some deterministic
penalties such as AIC or BIC.
These penalties are derived in some specific settings [e.g.,
a Gaussian noise is assumed by \citet{BiMa06}] and remain
unjustified and even sometimes misleading in more general settings.

Conversely, cross-validation (CV) is a \emph{resampling} procedure
based on a \emph{universal heuristics} which makes it applicable in a
wide range of settings.
CV procedures have been first studied in a
regression context by \citeauthor{Ston74} (\citeyear{Ston74,Ston77}) for
the leave-one-out
(Loo) and \citeauthor{Geis74} (\citeyear{Geis74,Geis75}) for
the $V$-fold cross-validation
(VFCV), and in the density estimation framework
by \citet{Rude82,Ston84}.
Since these procedures can be computationally demanding or even
intractable, \citet{Rude82,Bowm84} derived closed-form formulas for the
Loo estimator of the risk of histograms or kernel estimators.
These results have been recently extended to the leave-$p$-out
cross-validation (Lpo) by \citet{CeRo08}.

Although CV procedures are extensively used in practice, only few
theoretical results exist on their performances, most of them being of
asymptotic nature.
For instance, in the regression framework, \citeauthor{Burm89}
(\citeyear{Burm89,Burm90})
proves Loo is asymptotically the best CV procedure in terms
of risk estimation.
Several papers are dedicated to show the equivalence between some CV
procedures and penalized criteria in terms of asymptotic optimality
properties: (i) \emph{efficiency} in \citet{Li87,Zhan93}, and (ii)
\emph{model selection consistency} in \citet{Sha:1997}, \citeauthor{Yang06}
(\citeyear{Yang06,Yang07}).
Let us notice that in the parametric setting, \citet{Yang07} proved
that efficiency and model selection consistency are contradictory
objectives that cannot be achieved simultaneously.
We refer interested readers to
\citet{Sha:1997} for an extensive review about asymptotic optimality
properties in terms of efficiency and model selection consistency of some
penalized criteria as well as CV procedures.

As for nonasymptotic results in the density framework, \citet{BiMa97}
have settled an oracle inequality that relies on a conjecture and
may be applied to Loo.
However, to the best of our
knowledge, no such result has already been proved for Lpo in the
density estimation framework.
Recently, in the regression setting, \citet{Arlo07} established oracle
inequalities for $V$-fold penalties, while \citet{ArCe09} have
carried out an extensive simulation study in the change-point
detection problem with heteroscedastic observations.


In the present paper, we derive closed-form expressions for the Lpo
risk estimator of the broad class of projection estimators
(Section~\ref{sec.Lpo.closed.formulas}).
Such closed-form expressions considerably improve upon $V$-FCV in terms
of (i) variability [\citet{CeRo08}], and (ii) computational
complexity (Section~\ref{closed.form.expressions}).
A second improvement allowed by these formulas is the deep new
understanding of the theoretical performance of CV in two respects:
first for risk estimation (Section~\ref{subsec.risk.estimation}), and
second for model selection (Section~\ref{sec.model.selection}).
For instance, it is proved that Loo is the best CV procedure for risk
estimation (Theorem~\ref{thmm.optimal.risk.estimator}), while the story
can be different for model selection (Corollary~\ref
{cor.poly.coll.regularity} and Theorems~\ref
{thmm.model.consistency.belong} and~\ref{thmm.model.consistency.not.belong}).

In Section~\ref{sec.model.selection}, two aspects of model selection
\textit{via} CV have been explored.
The \emph{estimation} point of view is described in Section~\ref
{subsec.model.selection.CV.estimation}. It is shown that Lpo is optimal
as long as $p/n=o(1)$ and $p$ is large enough to balance the influence
of the model collection structure.
This phenomenon is supported by simulation experiments detailed in
Section~\ref{subsubsec.simulations}.
Finally, Section~\ref{subsec.identification} deals with the \emph
{identification} point of view.
CV is proved to be model selection consistent in various settings where
the choice of $p$ is related to the convergence rate (parametric and
nonparametric) of the best estimator one tries to recover.
Simulation results illustrate these different behaviors in Section~\ref
{subsubsec.Simuls.Identification}.
Finally, a discussion is provided in Section~\ref{sec.discussion} to
give some guidelines toward a better understanding of CV procedures.
The main proofs have been postponed to the \hyperref[app]{Appendix}.
For reasons owing to space constraints, more technical ones are
provided in the supplementary material [\citet{Celi:2014:supp}].

\section{Cross-validation and risk estimation}\label{sec.Lpo.closed.formulas}

\subsection{Statistical framework}\label{subsec.statistical.framework}

\subsubsection{Notation}
Throughout the paper, $X_1,\ldots,X_n\in[0,1]$ are independent and
identically distributed (\iid) random variables drawn from a probability
distribution $P$ of density $s\in L^2([0,1])$ with respect to
Lebesgue's measure on $[0,1]$, and $X_{1,n}=( X_1,\ldots,X_n )$.

Let $\mathcal{S}^*$ denote the set of measurable functions on
$[0,1]$. The distance between $\bayes$ and any $u\in\mathcal{S}^*$ is
measured by the quadratic \textit{loss} denoted by
\[
\ell\dvtx (s,u)\mapsto\ell( s,u ):=\| s-u \|^2 = \int
_{[0,1]} \bigl[ s(t)-u(t) \bigr]^2 \,dt.
\]
It is related to the \emph{contrast} function
%
\begin{equation}
\label{def.qudratic.contrast}\qquad \g\dvtx (u,x)\mapsto\g(u;x):=\| u \|^2-2u(x)
\qquad \mbox{with } \ell( s,u )=\Pg(u)-\Pg(s),
\end{equation}
where $\Pg(u)=P( \g(u;\cdot) )$ and $Pf:=\E[ f( X_1 ) ]$ for every
$f\in\mathcal{S}^*$.
The performance of an estimator $\sh=\sh(X_1,\ldots,X_n)$ of $\bayes$
is assessed by the \emph{quadratic risk}
\[
R_n(\sh):=\E\bigl[ \ell( s,\sh ) \bigr]=\E\bigl[ \| s-\sh
\|^2 \bigr].
\]

Estimating $\Pg(u)$ is made through the \emph{empirical
contrast} defined by
%
\begin{equation}
\label{exp.empirical.contrast} P_n\gamma(u):=\frac{1}{n}\sum
_{i=1}^n\g( u;X_i )\qquad\mbox {where }
P_n=1/n\sum_{i=1}^n
\delta_{X_i}
\end{equation}
denotes the empirical measure and $P_nf:=1/n\sum_{i=1}^n
f(X_i)$ for every $f\in\S^*$.

For a collection of models $\{ S_m \}_{m\in\mathcal{M}_n}$ indexed by
a countable
set $\M_n$, the \emph{empirical contrast minimizer} is defined by
%
\begin{equation}
\label{eq.empirical.contrast.minimizer.general} \widehat{s}_m:=\argmin_{u\in S_{m}}P_n
\gamma(u).
\end{equation}
It results a collection $\{ \widehat{s}_m \}_{m\in\mathcal{M}_n}$ of
estimators of
$\bayes
$ depending on the choice of models $S_m$s.
Instances of such models and estimators are described in Section~\ref
{subsubsec.projection.estimators}.

\subsubsection{Projection estimators} \label{subsubsec.projection.estimators}

Let $\Lambda_n$ be a set of countable indices and
$\{ \varphi_{\lambda} \}_{\lambda\in\Lambda_n}$ a family of
vectors in
$L^2([0,1])$ such that for every $m\in\mathcal{M}_n$,
$\{ \varphi_{\lambda} \}_{\lambda\in\Lambda(m)}$ denotes an
orthonormal family of
$L^2([0,1])$ with $\Lambda(m)\subset{\Lambda_n}$.
For every $m\in\mathcal{M}_n$, $S_{m}$ denotes the linear space
spanned by
$\{ \varphi_{\lambda} \}_{\lambda\in\Lambda(m)}$,
$D_m=\operatorname{dim}( S_m )$, and $s_m$ is the \emph{orthogonal
projection} of $\bayes$ onto $S_{m}$
\[
s_m:=\argmin_{u\in S_{m}}\Pg(u)=\sum
_{\lambda\in\Lambda
(m)}P\varphi_{\lambda} \varphi_{\lambda}\qquad
\mbox{with } P\varphi_{\lambda}=\E\bigl[ \varphi _{\lambda}(X) \bigr].
\]

\begin{defn}
An estimator $\sh\in L^2([0,1])$ is a projection estimator if there
exists a family $\{ \varphi_{\lambda} \}_{\lambda\in\Lambda}$ of
orthonormal
vectors of $L^2([0,1])$ such that
\[
\sh= \sum_{\lambda\in\Lambda}\theta_{\lambda}
\varphi_{\lambda
} \qquad\mbox{with } \theta_{\lambda}=\frac{1}{n}\sum
_{i=1}^n H_{\lambda}(X_i),
\]
where $\{ H_{\lambda}(\cdot) \}_{\lambda\in\Lambda}$ depends on the
family $\{ \varphi_{\lambda} \}_{\lambda\in\Lambda}$.
\end{defn}
It is straightforward to check that the empirical contrast minimizer
over $S_m=\operatorname{Span}( \varphi_{\lambda}, \lambda\in\Lambda(m)
)$, defined
by equation~\eqref{eq.empirical.contrast.minimizer.general}, is a
projection estimator since
%
\begin{equation}
\label{eq.projection.estimator.general} \widehat{s}_m=\sum_{\lambda\in\Lambda(m)}P_n
\varphi_{\lambda
}\varphi_{\lambda} \qquad\mbox {with } P_n
\varphi_{\lambda}=\frac{1}{n}\sum_{i=1}^n
\varphi_{\lambda
}(X_i).
\end{equation}
Here are a few examples of projection estimators [see \citet{DeLo93}]:
\begin{itemize}
\item\emph{Histograms}:
For every $m\in\mathcal{M}_n$, let $\{ I_{\lambda} \}_{\lambda\in
\Lambda(m)}$
be a partition of $[0,1]$ in $D_m = \card(\Lambda(m))$
intervals.
Set $\varphi_{\lambda}=\1_{I_{\lambda}}/\sqrt{\vert I_{\lambda}
\vert}$ for every
$\lambda\in\Lambda(m)$, with $\vert I_{\lambda} \vert$ the
Lebesgue measure
of $I_{\lambda}$, and $\1_{I_{\lambda}}(x)=1$ if $x\in I_{\lambda}$
and 0 otherwise. Then
%
\begin{equation}
\label{eq.histogram.estimator} \widehat{s}_m=\sum_{\lambda\in\Lambda(m)}P_n
\1_{I_{\lambda}} \frac{\1_{I_{\lambda}
}}{\vert I_{\lambda} \vert}.
\end{equation}

\item\emph{Trigonometric polynomials}:
For every $\lambda\in\Z$, let $\varphi_{\lambda}\dvtx t\mapsto
\varphi_{\lambda}
(t)=e^{2\pi
i\lambda t}$. Then for
any finite $\Lambda(m)\subset\mathbb{Z}$,
%
\begin{equation}
\label{eq.trigonometric.estimator} \widehat{s}_m(t)=\sum_{\lambda\in\Lambda(m)}P_n
\varphi_{\lambda
}e^{2\pi
i\lambda t}\qquad \forall t\in[0,1]
\end{equation}
is a trigonometric polynomial.

\item\emph{Wavelet basis}:
Let $\{ \varphi_{\lambda} \}_{\lambda\in\Lambda_n}$ be an
orthonormal basis of
$L^2([0,1])$ made of compact supported wavelets, where
$\Lambda_n=\{ (j,k)\mid j\in\mathbb{N}^* \mbox{ and } 1\leq
k\leq2^j \}$.
Then for every subset $\Lambda(m)$ of $\Lambda_n$,
%
\begin{equation}
\label{eq.wavelet.estimator} \widehat{s}_m=\sum_{\lambda\in\Lambda(m)}
P_n\varphi_{\lambda
}\varphi_{\lambda}.
\end{equation}
\end{itemize}
Some of these estimators can take negative values.
A possible solution is truncating and normalizing the preliminary
projection estimator
\[
\widetilde{s}_m = \widehat{s}_m\1_{\widehat{s}_m\geq0}
\biggl( \int_{[0,1]} \1_{\widehat{s}_m \geq0}(t) \widehat{s}_m(t)
\,dt \biggr)^{-1}.
\]
However, the closed-form expressions provided in Section~\ref
{closed.form.expressions} are not available for these truncated and
normalized estimators.

\subsection{Leave-$p$-out cross-validation}

In the literature, several cross-validation (CV) procedures have been
successively introduced to overcome the defects of already existing ones.
Let us describe the main CV procedures with some emphasis to
computational aspects.

\subsubsection{Cross-validation}\label{subsubsec.CV}

For $1\leq p\leq n-1$, let us define $\mathcal{E}_p=\{ e \subset \{
1,\ldots,n \},\break  \card(e)=p \}$ and for $e\in\mathcal{E}_p$, set
$X^e =
\{ X_i, i\in e \}$ (test set) and $X^{(e)} = \{ X_i, i\in\{ 1,\ldots,n \}\setminus e \}$ (training set).
Let also $P_n^{ e}:=1/p \sum_{i\in e} \delta_{X_i}$ and $P_n
^{ (e)}:=1/(n-p) \sum_{i\in(e)} \delta_{X_i}$ denote the
empirical measures, respectively, defined from the test set $X^e$ and the
training set $X^{(e)}$.

\paragraph*{Hold-out}

\emph{Simple validation} also
called \emph{Hold-out} was introduced in the early 1930s [\citet{Lars31}].
For every $1\leq p\leq n-1$, it consists in randomly splitting
observations into a training set $X^{(e)}$ of cardinality $n-p$ and a
test set $X^e$ of cardinality $p$. Random data splitting is only made
once and introduces additional variability.
For every $e\in\mathcal{E}_p$ (randomly chosen), the hold-out estimator
of $R_n(\sh)$ is
%
\begin{equation}
\label{def.hold.out.estimator} \Rh_{\mathrm{Ho},p}(\sh):= P_n^{ e}\g\bigl(
\sh\bigl( X^{(e)} \bigr) \bigr) = \frac{1}{p}\sum
_{i\in
e}\g\bigl( \sh\bigl(X^{(e)} \bigr);
X_i \bigr).
\end{equation}
Hold-out has been studied, for instance, by \citet{BaBL02,BlMa06} in
classification and by \citet{LuNo99,Wegk03} in regression.

\paragraph*{Leave-$p$-out}

Unlike equation~\eqref{def.hold.out.estimator} where a single split $e$
of the data is randomly chosen, which introduces additional unwanted
variability, \textit{leave-p-out} (Lpo) considers all the ${n\choose
p}= \card( \mathcal{E}_p )$ splits.
The Lpo estimator of $R_n(\sh)$ is defined by
%
\begin{equation}
\label{def.Lpo.estimator} \Rh_{p}(\sh)=\pmatrix{n\cr p}^{-1}\sum
_{e\in\mathcal{E}_p} P_n^{ e}\g\bigl( \sh\bigl(
X^{(e)} \bigr) \bigr).
\end{equation}
For instance, it has been studied by \citet{Shao93}, \citet{Zhan93},
and \citet{ArCe09} in the regression framework.
With $p=1$, Lpo reduces to the celebrated \emph{leave-one-out} (Loo)
cross-validation introduced by \citet{MoTu68} and further studied by
\citet{Ston74}.
Note that computing the Lpo estimator requires a computational
complexity of order ${n\choose p}$ times that of computing $\sh$, which
becomes intractable as $n$ grows.

\paragraph*{$V$-fold cross-validation}
To overcome the high computational burden of Lpo [equation~\eqref
{def.Lpo.estimator}], \citeauthor{Geis74} (\citeyear{Geis74,Geis75}) introduced the \emph{$V$-fold
cross-validation} ($V$-FCV).
Instead of considering all the ${n\choose p}$ possible splits, one
(randomly or not) chooses a partition of $X_1,\ldots,X_n$ into $V$
subsets $X^{e_1},\ldots,X^{e_V}$ of approximately equal size
$p=n/V=\card(e_i)$, $i=1,\ldots,V$.
Every $X^{e_i}$, $i=1,\ldots,V$ is successively used as a test set
leading to the $V$-fold risk estimator of $R_n(\sh)$
%
\begin{equation}
\label{def.VFCV.estimator} \Rh_{V\mathrm{\mbox{-}FCV}}(\sh)=\frac{1}{V}\sum
_{v=1}^V P_n^{ e_v}\g \bigl( \sh
\bigl( X^{(e_v)} \bigr) \bigr).
\end{equation}
$V$-FCV has been studied in the regression framework by
\citeauthor{Burm89}
(\citeyear
{Burm89,Burm90}) who suggests a correction to remove
its bias.

\subsubsection{Lpo versus $V$-FCV}\label{subsubsec.Lpo.VFCV}

As explained in Section~\ref{subsubsec.CV}, the Lpo computational
complexity is roughly ${n\choose p}$\vadjust{\goodbreak} times that of computing $\sh$,
which can be highly time-consuming.
Several surrogates of Lpo have been proposed such as $V$-FCV and the
\emph
{repeated learning-testing cross-validation} [\citet{BFOS84,Burm89,Zhan93}].
Unlike Lpo (and even Loo when $p=1$), $V$-FCV involves only $V$ such
computations, which is less demanding as long as $V\ll n$. Note that
usual values for $V$ are 3, 5, and 10 (except $V=n$ where $V$-FCV and Loo
coincide).\looseness=1

However, $V$-FCV relies on a preliminary (possibly random) partitioning
of $X_1,\ldots,X_n$ into $V$ subsets.
This preliminary partitioning induces some additional variability which
could be misleading.
For instance, \citet{CeRo08} have theoretically quantified the amount
of additional variability induced by $V$-FCV with respect to Lpo.

\subsection{Closed-form expressions for the Lpo risk estimator}
\label{closed.form.expressions}

%
Closed-form formulas for the Lpo estimator are proved in the present
section, which makes Lpo fully effective in practice and better than $V$-FCV.
Such formulas also enable a more accurate theoretical analysis of CV
procedures both in terms of risk estimation (Section~\ref
{subsec.risk.estimation}) and model selection (Section~\ref
{sec.model.selection}).

With the notation introduced at the beginning of Section~\ref
{subsubsec.CV}, let us consider projection estimators $\widehat{s}_m$
defined by
equation~\eqref{eq.projection.estimator.general}.
Closed-form formulas for the Lpo risk estimator are derived exploiting
the ``linearity'' of projection estimators.
Sums over $\mathcal{E}_p$ (which cannot be computed in general) then
reduce to binomial coefficients.
In the present section, proofs have been deferred to Appendix~A
[Supplementary material \citet{Celi:2014:supp}].
Recalling the expression of the contrast $\g(\cdot;\cdot)$
[equation~\eqref{def.qudratic.contrast}], one has to compute both
quadratic and linear terms.
%
\begin{lem}\label{lem.combi}
For every $m\in\mathcal{M}_n$, let $\widehat{s}_m=\widehat
{s}_m(X_{1,n})$ denote a projection
estimator defined by equation~\eqref{eq.projection.estimator.general} and
set $X^e = \{ X_i, i\in e \}$ for every $e\in\mathcal{E}_p$.
Then for every $p\in\{ 1,\ldots,n-1 \}$,
\begin{eqnarray*}
\sum_{e\in\mathcal{E}_p}\bigl\| \widehat{s}_m
\bigl(X^{(e)}\bigr) \bigr\|^2 & = &\frac{1}{(n-p)^2} \sum
_{\lambda\in\Lambda(m)} \Biggl[\pmatrix {n-1\cr p}\sum_{k=1}^n
\varphi_{\lambda}^2(X_k) \\
&&\hspace*{79pt}{}+\pmatrix{n-2\cr p}\sum
_{k\neq\ell} \varphi_{\lambda}(X_k)
\varphi_{\lambda
}(X_{\ell}) \Biggr],
\\
\sum_{e\in\mathcal{E}_p}\sum
_{i\in e}\sh\bigl(X^{(e)}\bigr) (X_i) & =&
\frac{1}{n-p} \sum_{\lambda\in\Lambda(m)} \pmatrix{n-2\cr p-1}\sum
_{i\neq
j}\varphi_{\lambda}(X_i)
\varphi_{\lambda}(X_j).
\end{eqnarray*}
\end{lem}
Lemma~\ref{lem.combi} enables to derive closed-form formulas for the Lpo
risk estimator, which makes Lpo procedure fully efficient in practice.
%
\begin{prop}\label{prop.global.proj.kernel}
For every $m\in\mathcal{M}_n$, let $\widehat{s}_m=\widehat
{s}_m(X_{1,n})$ denote a projection
estimator defined by equation~\eqref{eq.projection.estimator.general}.
Then for every $p\in\{ 1,\ldots,n-1 \}$,
%
\begin{eqnarray}\qquad
\label{exp.glob.proj.density} \Rh_p(m) &=& \Rh_p(\widehat{s}_m)
 \nonumber
 \\[-8pt]
 \\[-8pt]
 \nonumber
 &=& \frac{1}{n(n-p)}\sum_{\lambda
\in
\Lambda(m)} \Biggl[ \sum
_{j=1}^n\varphi_{\lambda}^2(X_j)-
\frac{n-p+1}{n-1}\sum_{j\neq k}\varphi_{\lambda}(X_j)
\varphi_{\lambda}(X_k) \Biggr].
\end{eqnarray}
\end{prop}
Proposition~\ref{prop.global.proj.kernel} enjoys a great interest.
First it applies to the broad family of projection estimators. Second,
it allows one to reduce the computation time from an exponential to a linear
complexity since computing \eqref{exp.glob.proj.density} is of order
$\mathcal{O}( n )$.
Note that in the more specific setting of histograms and kernel
estimators, such closed-form formulas have been derived by \citet{CeRo08}.

Let us now specify the Lpo estimator expressions for the three examples
of projection estimators given in Section~\ref
{subsubsec.projection.estimators}.
%
\begin{cor}[(Histograms)]
For $\widehat{s}_m$ given by equation~\eqref{eq.histogram.estimator}
and for
$p\in\{ 1,\ldots,n-1 \}$,
\[
\Rh_p(m)  = \frac{1}{(n-1)(n-p)}\sum
_{\lambda=1}^{D_m}\frac
{1}{\vert I_{\lambda} \vert} \biggl[ ( 2n-p )
\frac{n_{\lambda}}{n}-n(n-p+1) \biggl( \frac{n_{\lambda}}{n} \biggr)^2
\biggr],
\]
where $n_{\lambda}=\card( \{ i\st X_i\in I_{\lambda} \} )$.
\end{cor}
%
\begin{cor}[(Trigonometric polynomials)]
For every $k\in\N$, let $\varphi_{\lambda}$ denote either $
t\mapsto\cos(2\pi k
t)$, if $\lambda= 2k$ or $t\mapsto\sin(2\pi k t)$, if $\lambda=
2k+1$.
Let us further assume $\Lambda(m)=\{ 0,\ldots,2K \}$ for $K\in\N^*$.
Then for every $p\in\{ 1,\ldots,n-1 \}$,
\[
\Rh_p(m)  = \alpha(n,p) - \beta(n,p) \sum
_{k=0}^K \Biggl[ \Biggl\{ \sum
_{j=1}^n\cos(2\pi k X_j) \Biggr
\}^2+\Biggl\{ \sum_{j=1}^n
\sin(2\pi k X_j) \Biggr\}^2 \Biggr],
\]
where $\alpha(n,p) = (p-2)(K+1)[ (n-1)(n-p) ]^{-1}$ and $\beta(n,p)
= ( n-p+1 )[ n(n-1)(n-p) ]^{-1}$.
\end{cor}
%
\begin{cor}[(Haar basis)]
Let us define $\varphi(\cdot) = \1_{[0,1]}(\cdot)$ and
$\varphi_{j,k}(\cdot)=2^{j/2}\varphi(2^j\cdot-k)$, where $j\in
\mathbb{N}$ and $0\leq k\leq2^j-1$, and assume
$\Lambda(m)\subset\{ (j,k)\mid j\in\mathbb{N}, 0\leq k\leq2^j-1 \}$
for every $m\in\mathcal{M}_n$.
Then,
\[
\Rh_p(m)=\frac{1}{(n-1)(n-p)}\sum_{(j,k)\in\Lambda(m)}2^j
\biggl[ (2n-p)\frac{n_{j,k}}{n}-n(n-p+1) \biggl( \frac{n_{j,k}}{n}
\biggr)^2 \biggr],
\]
where $n_{j,k}=\card(  \{i\mid X_i\in[ k/2^j, (k+1)/2^j]
\} )$.
\end{cor}

\subsection{Risk estimation: Leave-one-out optimality}\label
{subsec.risk.estimation}

From the general\break formula~\eqref{exp.glob.proj.density}, one derives
closed-form expressions for the expectation and variance of the Lpo risk.
These expressions allow to analyze the theoretical behavior of CV in
terms of risk estimation and model selection (see Section~\ref
{sec.model.selection}).
In the present section, we prove the optimality of Loo for estimating
the risk of any projection estimator (Theorem~\ref
{thmm.optimal.risk.estimator}).

\begin{prop} \label{prop.expect.variance.proj.density}
For every $m\in\mathcal{M}_n$, let $\widehat{s}_m=\widehat
{s}_m(X_{1,n})$ denote a projection
estimator defined by equation~\eqref{eq.projection.estimator.general}.
Then for every $1\leq p\leq n-1$,
\[
\E\bigl[ \Rh_p(m) \bigr]  =\frac{1}{n-p}\sum
_{\lambda\in
\Lambda(m)}\bigl[ \E\varphi_{\lambda}^2(X)- \bigl(
\E\varphi_{\lambda}(X) \bigr)^2 \bigr]-\sum
_{\lambda\in\Lambda(m)}\bigl( \E \varphi_{\lambda}(X) \bigr)^2
\]
and
%
\begin{equation}
\label{exp.variance.general} \Var\bigl[ \Rh_p(m) \bigr]  = \frac{1}{(n-1)^2}
\biggl[ a_n + \frac {b_n}{(n-p)} + \frac{c_n}{(n-p)^2} \biggr],
\end{equation}
where
$ a_n = \Var[ \sum_{\lambda\in\Lambda(m)} ( n(P_n \varphi
_{\lambda})^2 -P_n\varphi_{\lambda}^2 ) ] $,
$c_n = \Var[ n \sum_{\lambda\in\Lambda(m)} ( P_n \varphi_{\lambda
}^2 - (P_n\varphi_{\lambda})^2 ) ] $, and $b_n = -2 \Cov[ \sum_{\lambda\in \Lambda (m)} ( n(P_n\varphi_{\lambda})^2 -P_n\varphi
_{\lambda}^2 ), \sum_{\lambda\in\Lambda(m)} n( P_n\varphi
_{\lambda}^2 - (P_n\varphi_{\lambda})^2 ) ] $.
\end{prop}
The proof is a straightforward application of Proposition~\ref
{prop.global.proj.kernel} and has been omitted.
Note that the above quantities do exist as long as $P\vert\varphi
_{\lambda} \vert^3<+\infty$ for any $\lambda\in\Lambda(m)$, which
holds true if
$\bayes$ is bounded for instance and
$\int\vert\varphi_{\lambda} \vert^3<+\infty$ ($\varphi_{\lambda
}$~continuous and compact supported,
e.g.).
In \eqref{exp.variance.general}, $a_n$, $b_n$ and $c_n$ do not depend
on~$p$.
Then knowing the behavior of the variance with respect to $p$ only
depends on the magnitude of $a_n$, $b_n$ and $c_n$, which is clarified
by Corollary~\ref{cor.variance.proj.density}.

Let us first focus on the bias $\B[ \Rh_p(m) ]:=\E\Rh
_p(m) -
\E[ \| \widehat{s}_m \|^2-2\int_{[0,1]}s \widehat{s}_m ]$ of the
Lpo estimator.
%
\begin{cor}[(Bias)]\label{cor.bias.proj.density}
For every $m\in\mathcal{M}_n$, let $\widehat{s}_m=\widehat
{s}_m(X_{1,n})$ denote a projection
estimator defined by equation~\eqref{eq.projection.estimator.general}.
Then for every $m\in\mathcal{M}_n$ and $1\leq p \leq n-1$,
\[
\B\bigl[ \Rh_p(m) \bigr]  = \frac{p}{n(n-p)}\sum
_{\lambda\in
\Lambda(m)}\Var\bigl[ \varphi_{\lambda}(X_1)
\bigr]\geq0.
\]
\end{cor}
The bias is nonnegative and increases with $p$, which means Loo ($p=1$)
has the smallest bias among CV procedures.
If $p=p_n$ satisfies $ p_n/n \mathop{\longrightarrow}\limits_{n\to+\infty} q \in[0,1)$,
then $\B[ \Rh_p(m) ] \mathop{\longrightarrow}\limits_{n\to+\infty} 0$, and Loo is
asymptotically unbiased.

Let us now describe the behavior of the variance with respect to $p$.
%
\begin{cor}[(Variance)]\label{cor.variance.proj.density}
With the same notation as Proposition~\ref
{prop.expect.variance.proj.density}, for every $m\in\mathcal{M}_n$
and $1\leq p
\leq n-1$,
\[
\Var\bigl[ \Rh_p(m) \bigr]  = \frac{n}{(n-1)^2} \biggl[ A +
\frac{ B }{n-p} + \frac{ C }{ (n-p)^2 } + O\biggl( \frac{1}{n} \biggr)
\biggr],
\]
where the big $O(\cdot)$ does not depend on $p$, but depends on $S_{m}$
and $P$, and
\begin{eqnarray*}
A & =& 4\Cov\biggl[ \sum_{\lambda}\varphi_{\lambda}(X_1)
\varphi_{\lambda
}(X_2), \sum_{\lambda }
\varphi_{\lambda}(X_1)\varphi_{\lambda
}(X_3)
\biggr] \geq0,
\\
B & =& 8\Cov\biggl[ \sum_{\lambda}\varphi_{\lambda}(X_1)
\varphi_{\lambda
}(X_2), \sum_{\lambda }
\varphi_{\lambda}(X_1)\varphi_{\lambda
}(X_3)
\biggr]\\
&&{} - 4 \Cov\biggl[ \sum_{\lambda}
\varphi_{\lambda}^2(X_1), \sum
_{\lambda
}\varphi_{\lambda} (X_1)
\varphi_{\lambda} (X_3) \biggr],
\\
C & =& 4\Cov\biggl[ \sum_{\lambda}\varphi_{\lambda}(X_1)
\varphi_{\lambda
}(X_2), \sum_{\lambda }
\varphi_{\lambda}(X_1)\varphi_{\lambda
}(X_3)
\biggr]\\
&&{} - 4 \Cov\biggl[ \sum_{\lambda}
\varphi_{\lambda}^2(X_1), \sum
_{\lambda
}\varphi_{\lambda} (X_1)
\varphi_{\lambda} (X_3) \biggr]
+ \Var\biggl[ \sum_{\lambda}\varphi_{\lambda}^2(X_1)
\biggr] \\
&\geq&0.
\end{eqnarray*}
\end{cor}
In the more specific case of histogram and kernel density estimators,
\citet{CeRo08} derived a similar (nonasymptotic) result for the variance.

The monotonicity of the variance with respect to $p$ depends on the
sign of $B$ since $x\mapsto f(x) = Ax^2+Bx+C$ has for derivative
$x\mapsto f'(x)=2Ax+B$ and $A\geq0$.
However, in full generality, the sign of $B$ is unknown.
The following proposition relates the monotonicity of $ p \mapsto\Var
[ \Rh_p(m) ] $ to this sign.
%
\begin{prop}\label{prop.variance.monotonicity}
Let us define $p_{0,n} = \argmin_{1\leq p\leq n-1} \Var[ \Rh _p(m)
]$ in equation~\eqref{exp.variance.general}. %
Then,
\[
\label{eqzerovariance} p_{0,n}  = n + \biggl( 1 - \frac{ \Cov[ \sum_{\lambda }\varphi_{\lambda
} ^2(X_1), \sum_{\lambda}\varphi_{\lambda}(X_1)\varphi_{\lambda
}(X_3) ] }{ 2\Cov[ \sum_{\lambda}\varphi_{\lambda}(X_1)\varphi
_{\lambda}(X_2), \sum_{\lambda}\varphi_{\lambda}(X_1)\varphi
_{\lambda} (X_3) ] } \biggr)
\bigl( 1+o(1) \bigr),
\]
where the little $o(\cdot)$ only depends on $S_m$ and $P$.
Furthermore, if
%
\begin{eqnarray}
\label{hypzerolocation} %
&&2\Cov\biggl[ \sum
_{\lambda}\varphi_{\lambda}(X_1)
\varphi_{\lambda
}(X_2), \sum_{\lambda }
\varphi_{\lambda} (X_1)\varphi_{\lambda
}(X_3)
\biggr]
\nonumber
\\[-8pt]
\\[-8pt]
\nonumber
&&\qquad\geq\Cov\biggl[ \sum_{\lambda}
\varphi_{\lambda}^2(X_1), \sum
_{\lambda}\varphi_{\lambda}(X_1)
\varphi_{\lambda}(X_3) \biggr],
\end{eqnarray}
%
%
$ p \in\{ 1,\ldots,n-1 \} \mapsto\Var[ \Rh_p(m) ] $ is increasing.
Otherwise,
$ p \mapsto\Var[ \Rh_p(m) ] $ is decreasing on $[1,p_{0,n}]$ and
increasing on $[p_{0,n}, n-1]$.
\end{prop}
Equation~\eqref{hypzerolocation} is related to the sign of $B$
(Corollary~\ref{cor.variance.proj.density}) and to the minimum location
value $p_{0,n}$.
If it holds true, then $p_{0,n}\notin\{ 2,\ldots,n-1 \}$, which means
Loo has the smallest variance among CV procedures.
In particular, let us notice \eqref{hypzerolocation} holds true with
any density estimated by regular histograms since $\sum_{\lambda
}\varphi_{\lambda}
^2(X_1)$ is then a constant and the covariance in the left-hand side is
a variance by independence of $X_1,X_2$, and $X_3$.
On the contrary, \eqref{hypzerolocation} is not fulfilled when using
histograms based on a partition such that $\P[ X\in I_{\lambda} ]=C<1/2$
for every $\lambda$, where $C$ denotes a constant.

We are now in position to provide the main result of this section,
which describes the behavior of $\Rh_p$ as a risk estimator in terms of
mean-square error (MSE).
%
\begin{thmm} \label{thmm.optimal.risk.estimator}
For every $m\in\mathcal{M}_n$, let us define the $\MSE$ of $\widehat
{s}_m$ by\break $\operatorname
{MSE}(m; p) = ( \B[ \Rh_p(m) ] )^2 + \Var[ \Rh_p(m) ]$,
for every $p\in\{ 1,\ldots,n-1 \}$.
\begin{longlist}[1.]
\item[1.] If \eqref{hypzerolocation} holds true,
then for every $m\in\mathcal{M}_n$, $p \mapsto\MSE(m; p)$ is
minimum for $p=1$.

\item[2.] Otherwise, for every $p=p_n\in\{ 1,\ldots,n-1 \}$ such that
$\lim
\sup_{n\to+\infty} p_n/n < 1$, then
\[
\MSE(m; p) = \frac{A}{n} + O\biggl( \frac{1}{n^2} \biggr)\qquad \mbox
{as } n\to+\infty,
\]
where $A$ is given in Corollary~\ref{cor.variance.proj.density} and the
big $O(\cdot)$ depends on $S_m$ and $P$.
\end{longlist}
\end{thmm}
If \eqref{hypzerolocation} holds true, Loo is the best CV procedure in
terms of MSE when estimating the risk of an estimator.
Otherwise as long as $\lim\sup_{n\to+\infty} p_n/n < 1$, choosing a
value of $p \neq1$ is useless since any value in $\{ 1,\ldots,n-1 \}$
asymptotically leads to the same performance in terms of MSE.
Therefore, since Loo has the smallest bias (Corollary~\ref
{cor.bias.proj.density}), \emph{Loo is optimal among CV procedures for
estimating the risk of an estimator}.
This result confirms what has been previously stated by \citet{Burm89}
in the regression framework.

\section{Optimal cross-validation for model selection}\label
{sec.model.selection}

From Section~\ref{subsec.risk.estimation}, Loo is proved to be the best
CV procedures in the context of risk estimation.
However, the best procedure for risk estimation is not necessarily the
best one for \emph{model selection}.
Although the empirical risk $P_n\gamma(\widehat{s}_m)$ \eqref
{eq.empirical.contrast.minimizer.general} is a reliable estimator of
$\E
[ P_n\gamma(\widehat{s}_m) ]$, using empirical risk minimization to
choose one
$\mh\in\mathcal{M}_n$ (without penalizing) would systematically lead
to overfitting.
The purpose of the present section is to study the performance of CV
for model selection with respect to the cardinality $p$ of the test set.

\subsection{Optimal cross-validation for estimation}\label
{subsec.model.selection.CV.estimation}

The performance of CV with respect to $p$ is first characterized
through a \emph{sharp oracle inequality} (Theorem~\ref
{thmm.poly.coll.regularity}).
A leading constant converging to 1 as $n\to+\infty$ is achieved for
some values of~$p$, highlighting the asymptotic optimality of
corresponding CV procedures.
From a theoretical point of view, Corollary~\ref
{cor.poly.coll.regularity} explores the link between (a proxy to) the
optimal $p$ and influential quantities related to the difficulty of the
estimation problem for finite sample size.
These results are further validated by simulation experiments
(Section~\ref{subsubsec.simulations}).

\subsubsection{Estimation point of view}

With the notation of Section~\ref{subsec.statistical.framework}, let us
consider a family of projection estimators $\{ \widehat{s}_m \}_{m\in
\mathcal{M}_n}$,
where $\mathcal{M}_n$ denotes an (at most countable) index set allowed
to depend
on $n$.
The best possible model, called \emph{the oracle model}, is denoted by
$S_{m^*}$, where
\begin{eqnarray*}
m^*&:=&\argmin_{m\in\mathcal{M}_n} \Pg(\widehat{s}_m) - \Pg(s) = \argmin
_{m\in\mathcal{M}_n} \| s - \widehat{s}_m \|^2
\\
& =& \argmin_{m\in\mathcal{M}_n} \Pg(\widehat{s}_m).
\end{eqnarray*}
Since $\Pg(\widehat{s}_m)$ has to be estimated, one uses CV (Lpo) to
choose a
candidate model for every $1\leq p \leq n-1$,
%
\begin{equation}
\label{exp.Lpo.index} \mh(p):=\argmin_{m\in\mathcal{M}_n} \Rh_p(m),
\end{equation}
and the final candidate model is denoted by $S_{\mh(p)}$.
The purpose is now to infer the properties of $\sh_{\mh(p)}$
with respect to $p\in\{ 1,\ldots,n-1 \}$ in terms of an oracle
inequality such as \eqref{exp.oracle.inequality.general}.

\subsubsection{Main oracle inequality}\label{subsubsec.oracle.inequality}
Let us introduce some notation and detail the main assumptions used
along the following sections.

\paragraph*{Square-integrable density}
%
\renewcommand{\theequation}{SqI}
\begin{equation}
\label{hyp.sq.integrable.density} s\in L^2\bigl([0,1]\bigr).
\end{equation}
Unlike \citet{Cast03}, for instance, it is not assumed that $s\geq
\rho
$ for a constant $\rho>0$.

\paragraph*{Polynomial collection}
There exists $a_{\M}\geq0$ such that
%
\renewcommand{\theequation}{Pol}
\begin{equation}
\label{hyp.poly.collection} \card( \mathcal{M}_n) \leq n^{a_{\mathcal{M}}}.
\end{equation}
In particular, this holds true if there exists $\alpha\geq0$ such that
$\card( \{ m\in\mathcal{M}_n, D_m=D \} ) \leq D^{\alpha}$, for every
$1\leq
D\leq n$.

\paragraph*{Model regularity}
%
\renewcommand{\theequation}{RegD}
\begin{equation}
\label{hyp.regular.basis.dimension} \exists\Phi>0,\qquad \sup_{m\in\mathcal{M}_n}
 \frac{ \| \phi_m \|
_{\infty} }{
D_m }
\leq\Phi\qquad \mbox{with } \phi_m = \sum_{\lambda\in
\Lambda(m)}
\varphi_{\lambda}^2. 
\end{equation}
It relates the regularity of the orthonormal basis (measured in terms
of sup-norm) to the dimension of the model.
For instance, using \eqref{eq.histogram.estimator}, \eqref
{hyp.regular.basis.dimension} requires $\vert I_{\lambda} \vert \geq
( \Phi
D_m)^{-1}$ for every $\lambda\in\Lambda(m)$. The length of intervals
$I_{\lambda}$ cannot be too different from one another to some extent.

\paragraph*{Maximal dimension}
%
\renewcommand{\theequation}{Dmax}
\begin{equation}
\label{hyp.dimension.maximal} \exists\Gamma>0,\qquad \sup_{m\in\mathcal{M}_n} D_m
\leq\Gamma\frac
{n}{(\log
n)^2}.
\end{equation}
In the sequel, $\Gamma=1$ is always considered to simplify expressions.
Note that proofs and conclusions remain unchanged with this particular choice.

\paragraph*{Estimation error and dimension}
%
\renewcommand{\theequation}{LoEx}
\begin{equation}
\label{hyp.expectation.lowerbound.dimension} \exists\xi>0,\qquad \inf_{m\in\mathcal{M}_n} \frac{\sqrt{n}\E( \| s_m
- \widehat{s}_m \| ) }{\sqrt{D_m}} \geq
\sqrt{\xi}. 
\end{equation}
This assumption makes the estimation error $\E( \| s_m - \widehat
{s}_m \|^2 )$ and $D_m$ comparable.
For instance, Lemma~B.3
[supplementary material \citet{Celi:2014:supp}] proves \eqref
{hyp.expectation.lowerbound.dimension} is fulfilled with H\"older
densities estimated by regular histograms defined by \eqref
{eq.histogram.estimator} such that $\vert\{ x \in[0,1] \mid
s(x)\geq \eta \} \vert\geq\ell$ for some $\eta\in(0,1)$, where
$0<\ell<1$ satisfies
\[
\ell> \Bigl( \inf_{m\in\mathcal{M}_n} D_m \Bigr)^{-1}.
\]
Note the latter inequality amounts to exclude too small models for
which the support of $s$ is included in one single interval $I_{\lambda}$.

\paragraph*{Richness of the collection}
There exist $m_0\in\M_n$ and $c_{\mathrm{rich}}\geq1$ such that,
%
\renewcommand{\theequation}{Rich}
\begin{equation}
\label{hyp.richness} \sqrt{n} \leq D_{m_0} \leq c_{\mathrm{rich}} \sqrt{n}.
\end{equation}
This requirement is rather mild since one can add such a model in our
collection.

\paragraph*{Approximation property}
There exist $c_{\ell},c_u>0$ and $\ell>u>0$ such that, for every
$m\in\M_n$,
%
\renewcommand{\theequation}{Bias}
\begin{equation}
\label{hyp.bias} c_{\ell} D_m^{-\ell} \leq\|
s-s_m \|^2 \leq c_u D_m^{-u}.
\end{equation}
This assumption quantifies the bias (approximation error) incurred by
model $S_{m}$ in estimating $\bayes$. It therefore relies on a smoothness
assumption on $\bayes$.
Such an upper bound is classical for $\alpha$-H\"olderian functions
with $\alpha\in(0,1]$ and regular histograms~\eqref
{eq.histogram.estimator}, for instance.
Note that \citet{Sto:1985} uses the same assumption (lower bound),
which is the \emph{finite sample counterpart} of the classical
assumption $\| s-s_m \|>0$ for every $m\in\mathcal{M}_n$ usually made
to prove
\emph{asymptotic optimality} for a model selection procedure
[see \citet{BiMa06}].

\paragraph*{Rate of convergence for the oracle model}
%
\renewcommand{\theequation}{OrSp}
\begin{equation}
\label{hyp.vitesse.risk.oracle} nR_n^* (\log n)^{-2} \mathop{\longrightarrow}\limits_{n\to
\infty}+\infty\qquad \mbox {with } R_n ^* = \inf_{m\in\mathcal{M}_n}
R_n(\widehat{s}_m). 
\end{equation}
The risk of the oracle model $R_n^*$ does not decrease to 0 faster than
$ (\log n)^2/n$.
In particular, this holds true for densities in $\mathcal{H}(L,\alpha)$
with $L>0$ and $\alpha\in(0,1]$ estimated by regular histograms [see
Section~B.5 in \citet{Celi:2014:supp}].

The performance of the Lpo estimator with respect to $p$ is described
by the following oracle inequality from which the CV optimality is
deduced for some values of $p$.
The proof is given in Appendix~\ref{app.proof.oracle}.
%
\begin{thmm}[(Optimal CV)]
\label{thmm.poly.coll.regularity}
Let $\bayes$ denote a density on $[0,1]$ such that \eqref
{hyp.sq.integrable.density} holds true, set $\{ S_m \}_{m\in\M_n}$ a
collection of models defined in Section~\ref
{subsubsec.projection.estimators}, and assume \eqref
{hyp.poly.collection}, \eqref{hyp.regular.basis.dimension}, \eqref
{hyp.dimension.maximal}, \eqref{hyp.richness}, \eqref
{hyp.expectation.lowerbound.dimension}, \eqref{hyp.bias} and \eqref
{hyp.vitesse.risk.oracle}.
Let $\mh=\mh(p)$ denote the model minimizing $\Rh_p(m)$ over $\M_n$ for
every $p\in\{ 1,\ldots,n-1 \}$.
Then there exist a sequence $(\delta_n)_{\N}$ such that $\delta_n\to
0$, and $n\delta_n\to+\infty$ as $n\to+\infty$, and an event
$\widetilde{\Omega}$ with $\P(\widetilde{\Omega}) \geq1- 6/n^2$ on
which, for large enough values of $n$,
\[
\| s-\sh_{\mh(p)} \|^2 \leq C_{n}(p) \inf
_{m\in\mathcal{M}_n} \bigl\{ \| s-\sh _{m} \|^2 \bigr\}\qquad
\mbox{with } C_{n}(p) = \frac{ T_{B}^+ \vee
T_{V}^+ }{ T_{B}^- \wedge T_{V}^- } \geq1,
\]
where
\begin{eqnarray*}
T_{B}^- & =& 1 - \delta_n K(n,p), \\
 T_{V}^- &=&
\frac
{1}{1-p/n}(1-\delta_n)[ 1 - 4 \delta_n ] - 2
\delta_n K(n,p) [ 3 - 4\delta_n ],
\\
T_{B}^+ & = &1 + \delta_n K(n,p), \\
 T_{V}^+ &=&
\frac{1}{1-p/n}(1+\delta_n)[ 1 + 4 \delta_n ] + 2\delta
_n K(n,p) [ 3 + 4\delta_n ],
\end{eqnarray*}
and $K(n,p) = 1 + \frac{2}{n-1} + \frac{p}{n-p} \frac{1}{n-1}$.
\end{thmm}
If $p/n\to0$ then $C_n(p) \to1$ as $n\to+\infty$, which leads to
\emph{efficient} (asymptotically optimal) model selection procedures.
In particular, this holds true for $p=1$ that is, \emph{Loo is
asymptotically optimal} since
\[
\frac{ \| s-\sh_{\mh(1)} \|^2 }{ \inf_{m\in\mathcal{M}_n} \{ \|
s-\sh _{m} \|^2 \} } \mathop{\longrightarrow}\limits_{n\to+\infty}^{\mathrm{a.s.}} 1.
\]
From the proof, it also arises that the slowly decreasing sequence $\delta
_n$ is related to the model collection structure. An increase of $\M_n$
makes the model selection problem more difficult and $\delta_n$ larger.

\begin{figure}

\includegraphics{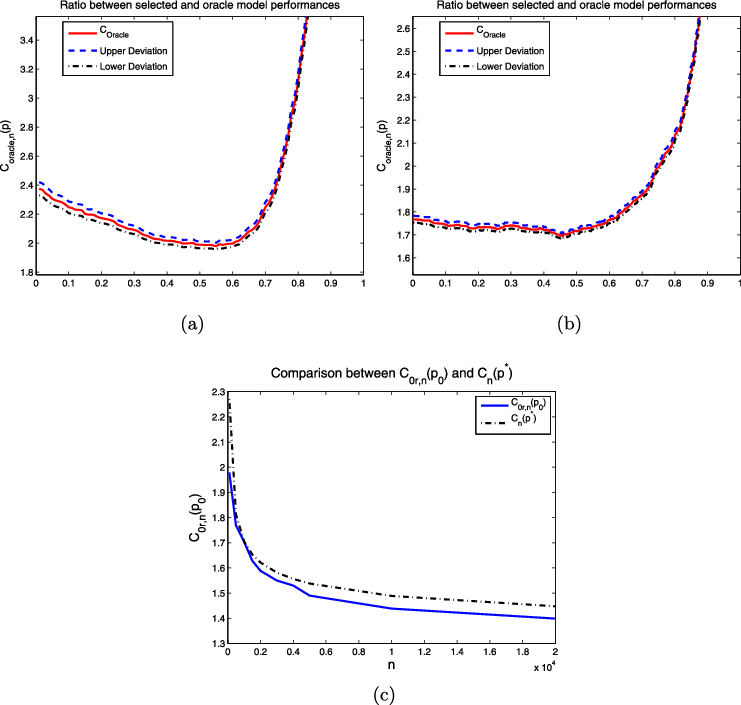}

\caption{Panels~\textup{(a)} and~\textup{(b)}: $p/n
\mapsto C_{\mathrm{or},n}(p)$ (plain red line) is plotted for $\Gamma=1$ [see
\protect\eqref{hyp.dimension.maximal}] and different values of $n$: \textup{(a)}
$n=100$, and \textup{(b)} $n=1000$.
$p/n \mapsto C^+_{\mathrm{oracle},n}(p)$ (blue dashed line) and $p/n \mapsto
C^-_{\mathrm{oracle},n}(p)$ (black dot-dashed line) have been plotted on the
same graph as well [see \protect\eqref{exp.upper.lower.deviation}].
Panel~\textup{(c)}: $n \mapsto C_{\mathrm{or},n}(p_0)$ (plain blue line) and $n \mapsto
C_{n}(p^*)$ (black dot-dashed line) are displayed.
$N=1000$ samples have been drawn from the mixture of Beta distributions
\protect\eqref{exp.mixture.beta.distribution}.}\label{fig.sample.size.Gamma1}
\end{figure}

While \emph{asymptotic optimality} is deduced from Theorem~\ref
{thmm.poly.coll.regularity} for any CV procedure as long as $p=o(n)$, it
is also desirable to analyze the performance of CV as $p$ depends on
the finite sample size. From a theoretical point of view, this will
provide the rate at which $p/n$ has to decrease to 0 to reach efficiency.
Based on Figure~\ref{fig.sample.size.Gamma1} [panel~(c)] where $C_n(p)$
appears as a reliable proxy to the optimal $C_{\mathrm{or},n}(p)$ [given by
equation~\eqref{exp.optimal.ratio}], we suggest to optimize
$C_n(p)$ with respect to $p$ to get a surrogate optimal $p$ depending
on influential parameters such as $n$ and $\delta_n$.
This strategy has been validated by simulation experiments of
Section~\ref{subsubsec.simulations}.
The following Corollary~\ref{cor.poly.coll.regularity} proves the best
(surrogate) $p/n$ slowly decreases to $0$.

\begin{figure}

\includegraphics{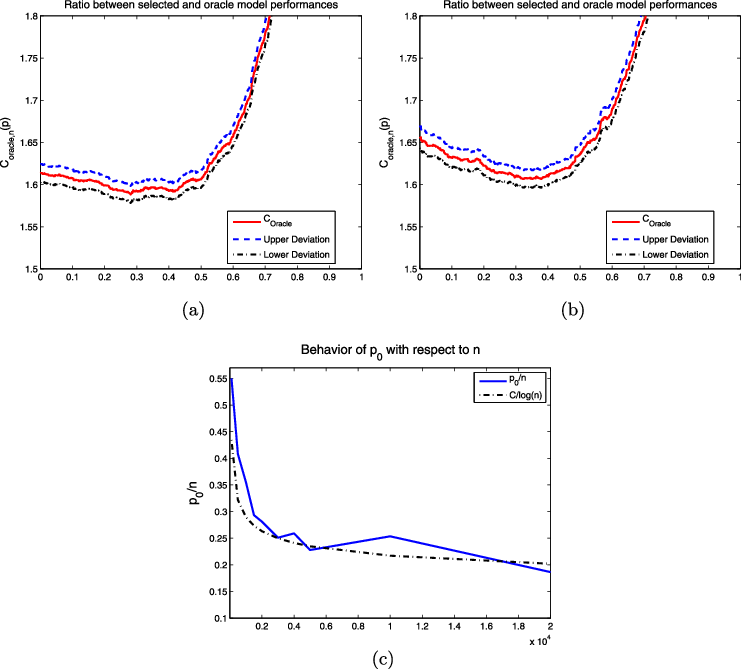}

\caption{For \textup{(a)} and \textup{(b)}, $p/n
\mapsto C_{\mathrm{or},n}(p)$ (plain red line) is plotted for $n=2000$ and
different values of $\Gamma$ [see \protect\eqref
{hyp.dimension.maximal}]: \textup{(a)}
$\Gamma=1$, \textup{(b)} $\Gamma=2$.
$p/n \mapsto C^+_{\mathrm{oracle},n}(p)$ (blue dashed line) and $p/n \mapsto
C^-_{\mathrm{oracle},n}(p)$ (black dot-dashed line) have been plotted on the
same graph as well [see \protect\eqref{exp.upper.lower.deviation}].
$N=1000$ samples have been drawn from the mixture of Beta distributions
\protect\eqref{exp.mixture.beta.distribution}.
For \textup{(c)}, $n \mapsto p_0/n$ (blue plain line) and $n\mapsto C/(\log n)$
(black dot-dashed line) are displayed, where $p_0$ denotes the
minimizer of $C_{\mathrm{or},n}(p)$ as a function of $p$ and $C$ is a constant.}
\label{fig.sample.size.Gamma2.popt}
\end{figure}

%
\begin{cor}[(Optimizing upper bound)]
\label{cor.poly.coll.regularity}
With the notation and assumptions of Theorem~\ref
{thmm.poly.coll.regularity}, the constant $C_n(p)$ is minimized over
$p\in\{ 1,\ldots,\break n-1 \}$
for
\[
0< \frac{p_n^*}{n} = 1 - \frac{1 - 5\delta_n +4\delta_n^2 -
({2}/{(n-1)}) (3\delta_n - 4\delta_n^2) + {\delta_n}/{(n-1)} }{1 +
2(1+{1}/{(n-1)}) (3\delta_n - 4\delta_n^2) - \delta_n(1+{1}/{(n-1)}) } < 1.
\]
Furthermore, the optimal ratio $p^*/n$ is slowly decreasing to 0 as $n$
tends to $+\infty$
%
\renewcommand{\theequation}{\arabic{equation}}
\setcounter{equation}{15}
\begin{equation}
\label{lim.p.optimal.rate} p^*_n \sim_{+\infty} 10 n\delta_n
\mathop{\longrightarrow}\limits_{n \to+\infty} +\infty.
\end{equation}
\end{cor}
The proof has ben deferred to Appendix~\ref{app.proof.oracle}.
Corollary~\ref{cor.poly.coll.regularity} describes the rate (up to
constant) at which $p=p_n^*$ has to grow with $n$ to achieve
finite-sample optimality.
%
%
In particular $p_n^*/n$ in \eqref{lim.p.optimal.rate} is related to
$\delta_n$ which is strongly connected to the structure of the model
collection as explained following Theorem~\ref
{thmm.poly.coll.regularity}.
A~more complex collection leads to a larger $\delta_n$ and then to a
larger optimal $p_n^*$. In other words, $p$ must be chosen large enough
to balance the overfitting induced by the structure of the model collection.
This phenomenon is observed in practice in the simulation experiments
of Section~\ref{subsubsec.simulations} (Figure~\ref
{fig.sample.size.Gamma2.popt}).

%

\subsubsection{Adaptivity in the minimax sense}

\emph{Adaptivity in the minimax sense} is a desirable property for
model selection procedures. It means the considered procedure
automatically adapts to the unknown smoothness of the target function
$s$ to estimate
[see \citet{BarBirMas:1999} for an extensive presentation].

Several adaptivity in the minimax sense results are provided in the
present section.
Deriving such results from oracle inequalities \eqref
{exp.oracle.inequality.general} is somewhat classical.
Here, the novelty is first that CV enjoys such a desirable property as
a model selection procedure, second that the leading constant $C_n(p)$
in Theorem~\ref{thmm.poly.coll.regularity} when converging to 1 as $n$
tends to $+\infty$ provides accurate results.

Let us start providing a general theorem from which any adaptivity
result will be immediate corollary. The proof is given in Appendix~\ref
{app.proof.oracle}.
%
\begin{thmm}\label{thmm.oracle.inequality.expectation}
Let $\bayes$ denote a density on $[0,1]$ such that \eqref
{hyp.sq.integrable.density} holds true, set $\{ S_m \}_{m\in\M_n}$ a
collection of models defined in Section~\ref
{subsubsec.projection.estimators}, and assume \eqref
{hyp.poly.collection}, \eqref{hyp.regular.basis.dimension}, \eqref
{hyp.dimension.maximal}, \eqref{hyp.richness}, \eqref
{hyp.expectation.lowerbound.dimension}, \eqref{hyp.bias}, and \eqref
{hyp.vitesse.risk.oracle}.
Let $\mh=\mh(p)$ denote the model minimizing $\Rh_p(m)$ over $\M_n$ for
every $p\in\{ 1,\ldots,n-1 \}$.
Then for every $1\leq p \leq n-1$,
%
\renewcommand{\theequation}{\arabic{equation}}
\setcounter{equation}{16}
\begin{eqnarray}
\label{ineq.oracle.inequality.expectations} &&\E\bigl[ \| \bayes- \sh_{\mh(p)} \|^2 \bigr]
\nonumber
\\[-8pt]
\\[-8pt]
\nonumber
&&\qquad\leq C_n(p) \E\Bigl[ \inf_{m\in
\mathcal{M}_n} \| \bayes-
\sh_{m} \|^2 \Bigr] + \bigl( \Phi+ \| s \|^2
\bigr) \frac{12}{n (\log n)^2} + \frac{6 c_u}{n^2},
\end{eqnarray}
where $C_{n}(p) = \frac{ T_{B}^+ \vee T_{V}^+ }{ T_{B}^- \wedge T_{V}^-
}$, with
\begin{eqnarray*}
T_{B}^- & =& 1 - \delta_n K(n,p),\\
  T_{V}^- &=&
\frac
{1}{1-p/n}(1-\delta_n)[ 1 - 4 \delta_n ] - 2
\delta_n K(n,p) [ 3 - 4\delta_n ],
\\
T_{B}^+ & =& 1 + \delta_n K(n,p),\\
 T_{V}^+ &=&
\frac{1}{1-p/n}(1+\delta_n)[ 1 + 4 \delta_n ] + 2\delta
_n K(n,p) [ 3 + 4\delta_n ],
\end{eqnarray*}
and $K(n,p) = 1 + \frac{2}{n-1} + \frac{p}{n-p} \frac{1}{n-1}$.
\end{thmm}
The last two terms in the right-hand side of \eqref
{ineq.oracle.inequality.expectations} are remainder terms by
Assumptions \eqref{hyp.regular.basis.dimension}, \eqref
{hyp.dimension.maximal}, and \eqref{hyp.bias}.

Applying Theorem~\ref{thmm.oracle.inequality.expectation} to the
collection of regular histograms defined by \eqref
{eq.histogram.estimator}, the following corollary settles an adaptivity
property with respect to H\"older balls [see \citet{DeLo93}].
%
\begin{cor}\label{cor.adaptivity.holder}
Let us consider the model collection of Section~\ref
{subsubsec.projection.estimators} made of piecewise constant functions
and the associated histograms defined by \eqref{eq.histogram.estimator}
such that, for every $m\in\mathcal{M}_n$ and $\lambda\in\Lambda
(m)$, $\vert I_{\lambda} \vert=D_m^{-1}$ (regular histograms).
Let us also assume \eqref{hyp.dimension.maximal} and \eqref
{hyp.expectation.lowerbound.dimension} hold true.

If the target density $\bayes$ belongs to the H\"older ball $\mathcal
{H}(L,\alpha)$ for some $L>0$ and $\alpha\in(0,1]$, then there exist
constants $0<K_{\alpha}^-\leq K_{\alpha}^+$ such that for every $ p =
o(n) $,
\begin{eqnarray*}
&&K_{\alpha}^- L^{{2}/{(2\alpha+1)}} n^{-{2\alpha}/{(2\alpha+1)}} \\
&&\qquad\leq\sup
_{s\in\mathcal{H}(L,\alpha)} \E\bigl[ \| s- \sh _{\mh (p)} \| ^2
\bigr] \\
&&\qquad\leq\bigl( 1 + o(1) \bigr) K_{\alpha}^+ L^{{2}/{(2\alpha+1)}}
n^{-{2\alpha}/{(2\alpha+1)}} + O\biggl( \frac{1}{n (\log n)^2} \biggr),
\end{eqnarray*}
$K_{\alpha}^-$ and $K_{\alpha}^+$ only depend on $\alpha$ (not on $n$
or $s$).

Furthermore, since this property holds for every $L>0$ and $\alpha\in
(0,1]$, then $\{ \sh_{\mh(p)} \}_{n\in\N^*}$ is adaptive in the minimax
sense with respect to $\{ \mathcal{H}(L,\alpha) \}_{L>0,\alpha\in
(0,1]}$ for every $ p = o(n) $.
\end{cor}
The proof has been deferred to Section~B.5 in \citet{Celi:2014:supp}.
The upper bound is tight since the rate $ n^{-{2\alpha}/{(2\alpha
+1)}}$ and the dependence on the radius $L^{{2}/{(2\alpha+1)}}$ are
the same as in the lower bound, which has been stated by \citet{IbKh81}.
The main contribution of this result is to prove $p=o(n)$ leads to adaptivity.
Note that similar results can also be proved for Besov balls $\mathcal
{B}_{\infty,2}^{\alpha}(L)$, with $\alpha,L>0$, for instance
[see \citet{DeLo93}], by using an appropriate collection of models such as
trigonometric polynomials defined by~\eqref{eq.trigonometric.estimator}.

\subsubsection{Simulation experiments}\label{subsubsec.simulations}

Results of simulation experiments are provided to check the conclusions
drawn (from theory) in Section~\ref{subsubsec.oracle.inequality}.
A mixture of Beta distributions
%
\renewcommand{\theequation}{\arabic{equation}}
\setcounter{equation}{17}
\begin{equation}
\label{exp.mixture.beta.distribution} \forall x\in[0,1]\qquad  s(x) = \frac{\beta(3,7; x) + \beta(10,5; x) }{2}
\end{equation}
has been used to generate samples of size
$n=100,500,1000,2000,3000,4000,\break 5000,10\mbox{,} 000,20\mbox{,} 000$.
Note that \eqref{exp.mixture.beta.distribution} defines a H\"older
density on $[0,1]$.
For each $n$, every $p \in\{ 1,\ldots,n-1 \}$ have been considered and
\eqref{hyp.dimension.maximal} is fulfilled with $\Gamma= 1$
(Figure~\ref{fig.sample.size.Gamma1}) and $\Gamma=2$ (Figure~\ref
{fig.sample.size.Gamma2.popt}).

The model collection we used is made of piecewise constant functions
described in Section~\ref{subsubsec.projection.estimators} leading to
regular histogram estimators defined by \eqref{eq.histogram.estimator}.
Only regular histograms with dimension $D_m\geq2$ are used so that
\eqref{hyp.expectation.lowerbound.dimension} holds true (Lemma~B.3).
For every $1\leq p \leq n-1$, $\mh(p)$ is defined by \eqref{exp.Lpo.index}.

Let us also introduce
%
\begin{equation}\qquad
\label{exp.optimal.ratio} C_{\mathrm{or},n}(p):=\E\biggl[ \frac{ \| s-\sh_{\mh(p)} \|^2 }{ \inf_{m\in
\mathcal{M}_n} \{ \| s-\sh_{m} \|^2 \} } \biggr]\quad
\mbox{and} \quad p_0:=\argmin_{1\leq p \leq n-1} C_{\mathrm{or},n}(p),
\end{equation}
which measures the average performance of $\sh_{\mh(p)}$ with respect
to that of $\sh_{m^*}$ (oracle estimator). The closer $C_{\mathrm{or},n}(p)$ to
1, the better $ \sh_{\mh(p)}$.
Minimizing $C_{\mathrm{or},n}(p)$ as a function of $p$ for various values of $n$
enables to check whether the conclusions drawn from minimizing $C_n(p)$
with respect to $p$ (Theorem~\ref{thmm.poly.coll.regularity} and
Corollary~\ref{cor.poly.coll.regularity}) hold true or not, that is
whether $C_n(p)$ is an accurate approximation to $C_{\mathrm{or},n}(p)$.
For each curve $p\mapsto C_{\mathrm{or},n}(p)$, a confidence band has been
displayed. It is delimited by $p\mapsto C_{\mathrm{or},n}^-(p)$ and $p\mapsto
C_{\mathrm{or},n}^+(p)$, respectively, defined by
%
\begin{equation}
\label{exp.upper.lower.deviation} C_{\mathrm{or},n}^-(p) = C_{\mathrm{or},n}(p) -
\frac{\widehat{\sigma}}{\sqrt{N}}\quad\mbox{and}\quad C_{\mathrm{or},n}^+(p) = C_{\mathrm{or},n}(p) +
\frac{\widehat{\sigma
}}{\sqrt{N}},
\end{equation}
where $\widehat{\sigma}$ denotes the empirical standard deviation.

%

First, from panels~(a) and~(b) of Figure~\ref{fig.sample.size.Gamma1},
$ C_{\mathrm{or},n}(p)$ (plain red lines) decreases pointwise as $n$ grows. This
is confirmed by panel~(c) of Figure~\ref{fig.sample.size.Gamma1} at the
particular value $p=p_0$ as $n$ grows. This is in accordance with
Theorem~\ref{thmm.poly.coll.regularity} and $C_n(p)\to0$ as $n$ increases.
Second, the optimization strategy at the basis of Corollary~\ref
{cor.poly.coll.regularity} is empirically validated
by panel~(c) of Figure~\ref{fig.sample.size.Gamma1} where
$C_{\mathrm{or},n}(p_0)$ and its proxy $C_{n}(p_0)$ remain very close to each
other. Furthermore, the optimal rate derived in Corollary~\ref
{cor.poly.coll.regularity} is supported up to constant by simulation
results displayed in panel~(c) of Figure~\ref
{fig.sample.size.Gamma2.popt} where $p_0/n$ is almost equal to the
predicted $\delta_n\approx C/(\log n)$ ($C>0$) from the proof of
Theorem~\ref{thmm.poly.coll.regularity}.

The conclusion of Corollary~\ref{cor.poly.coll.regularity} about the
dependence of the optimal $p_0$ on the complexity of the model
collection (through $\delta_n$) is also illustrated by panels~(a)
and~(b) in Figure~\ref{fig.sample.size.Gamma2.popt} where $\Gamma$
\eqref{hyp.dimension.maximal}, respectively, equals $1$ and $2$.
As $\Gamma$ grows the model collection becomes more complex, leading to
a worse performance and a larger $p_0$ in panel~(b). The need for a
larger $p_0$ is all the more strong as the curve in panel~(b) is less
flat than in panel~(a), indicating the problem becomes more difficult
as $\Gamma$ increases and any misspecification of $p_0$ leads to a
stronger loss in accuracy.
One concludes the more complex the model collection, the larger the
optimal $p$.

Note that this conclusion does not apply to Loo ($p=1$) [see \eqref
{lim.p.optimal.rate}], suggesting Loo may be suboptimal for finite
sample size. This is supported by Figure~\ref{fig.sample.size.Gamma1}
[panels~(a) and (b)] and Figure~\ref{fig.sample.size.Gamma2.popt}
[panels~(a) and (b)] where the minimum of each curve is not reached at $p=1$.

%
%
%

\subsection{Optimal cross-validation for identification}
\label{subsec.identification}


With the notation of Section~\ref{subsec.statistical.framework}, $\{
\widehat{s}_m \}_{m\in\mathcal{M}_n}$ denotes a collection of
projection estimators
(Section~\ref{subsubsec.projection.estimators}) which is allowed to
depend on $n$.
The purpose is now to recover the best model denoted by $S_{\bar m}$
and defined by
%
\begin{equation}
\label{expr.true.model.identification} \bar m:=\argmin_{m\in\mathcal{M}_n} \E\bigl[ \| s -
\widehat{s}_m \|^2 \bigr],
\end{equation}
where $\bar m$ is a deterministic quantity unlike $m^*$ from
Section~\ref{subsec.model.selection.CV.estimation}.
Since this goal cannot be reached if other models can perform as well
as $S_{\bar m}$ (even asymptotically), one also requires there exist
$\mu>0$ and $n_0\in\N^*$ such that for every integer $n>n_0$,
%
\renewcommand{\theequation}{BeMo}
\begin{equation}
\label{hyp.best.model.identification} (1+\mu) \E\bigl[ \| s - \sh_{\bar m} \|^2
\bigr]\leq\inf_{m\in\M
_n\setminus\{ \bar m \}} \E\bigl[ \| s - \widehat{s}_m
\|^2 \bigr]. 
\end{equation}
A similar assumption (in probability rather than in expectation) has
been made by \citet{Yang07}.
Let us further assume the collection $\{ S_m \}_{m\in\M_n}$ can be
split into:
\begin{itemize}
\item\emph{parametric models} indexed by $\M_{n,P}$ for which there
exist constants $\pi,\tau>0$ (independent of $n$) such that
%
\renewcommand{\theequation}{\arabic{equation}}
\setcounter{equation}{21}
\begin{equation}
\label{hyp.PNP.parametric} \sup_{m\in\mathcal{M}_{n,P} } \bigl\{ n\E\bigl[ \|
s_m - \sh _{ m} \|^2 \bigr] \bigr\} \leq\pi\quad
\mbox{and}\quad \inf_{m\in\mathcal
{M}_{n,P}, s\notin S_m } \bigl\{ \| s - s_m
\|^2 \bigr\} \geq\tau.
\end{equation}

\item\emph{nonparametric models} indexed by $\M_{n,NP}$ such that
%
\begin{equation}
\label{hyp.PNP.nonparametric} n ( \log n )^{-2} \inf_{m\in\mathcal{M}_{n,NP}} \E\bigl[
\| s_m - \sh _{m} \|^2 \bigr] \mathop{\longrightarrow}\limits_{n
\to+\infty} +\infty.
\end{equation}
Then
%
\renewcommand{\theequation}{P-NP}
\begin{equation}
\label{hyp.param.nonparam} 
\{ S_m \}_{m\in\M_n} = \{
S_m \}_{m\in\M_{n,P}} \cup\{ S_m \} _{m\in\M_{n,NP}}.
\end{equation}
\end{itemize}
Parametric models are models with convergence rate of order $1/n$.
Since $\E[ \| s - \sh_{ m} \|^2 ] \approx\| s-s_m \|^2 +
C\cdot D_m/n$, allowing $D_m$ to depend on $n$ makes the rate of the
corresponding model slower than $1/n$ (nonparametric model).
Consistently with this remark, \eqref{hyp.PNP.parametric} requires the
largest dimension over parametric models is bounded by a constant
independent of $n$, and that the bias of parametric models such that
$s\notin S_m$ cannot decrease with $n$ toward 0. Otherwise, such a
model would be nonparametric.
Conversely, \eqref{hyp.PNP.nonparametric} only requires that the dimension
of nonparametric models must be larger than $(\log n)^2$. In
particular, this does not prevent nonparametric models from containing
$s$ or having their bias decreasing to 0 as $n$ grows.

\subsubsection{Main results}

Depending on whether $s$ belongs or not to $\bigcup_{m\in\M_n} S_{m}$, the
two following results prove model selection consistency for CV.
Their main contribution is to relate the cardinality $p$ of the test
set to the rate of convergence of $\sh_{\bar m}$ and the model
collection complexity.
Note that, in addition, the model consistency property is settled with a
collection of models allowed to grow with $n$, which contrasts with
earlier results [see, e.g., \citet{Yang07}].

Let us start with the setting where $s$ belongs to $\bigcup_{m\in\M_n}
S_{m}
$, which implies the best estimator $\sh_{\bar m}$ achieves the
parametric rate $1/n$.
%
\begin{thmm}[(Model consistency with $s \in\bigcup_m S_{m}$)] 
\label{thmm.model.consistency.belong}
Let $\bigcup_{m\in\M_n} S_{m}$ denote a collection of models
satisfying~\eqref{hyp.poly.collection} and~\eqref{hyp.param.nonparam},
$\bar m \in\mathcal{M}_n$ given by~\eqref
{expr.true.model.identification} be such
that~\eqref{hyp.best.model.identification} holds true, and
assume~\eqref
{hyp.sq.integrable.density},~\eqref
{hyp.regular.basis.dimension},~\eqref
{hyp.dimension.maximal}, and~\eqref{hyp.expectation.lowerbound.dimension}.
For every $1\leq p\leq n-1$, let us also define $\mh= \mh(p) =\break
\argmin
_{m\in\M_n} \Rh_p(m)$.
If the target $s \in\bigcup_{m\in\M_n} S_{m}$, then every $1\leq p=p_n
\leq
n-1$ such that
%
\renewcommand{\theequation}{\arabic{equation}}
\setcounter{equation}{23}
\begin{equation}
\label{conditions.p.belong} \log(n) \biggl( 1- \frac{p}{n} \biggr)\mathop{\longrightarrow}\limits_{n\to+
\infty} 0 \quad\mbox{and}\quad n\biggl( 1- \frac{p}{n} \biggr) \mathop{\longrightarrow}\limits_{n
\to+\infty} +\infty,
\end{equation}
leads to
\[
\P[ \mh= \bar m ] \mathop{\longrightarrow}\limits_{n\to+\infty} 1.
\]
\end{thmm}
The proof has been deferred to Appendix~\ref
{subsec.main.proofs.identification}.
When $s$ belongs to $\bigcup_{m\in\M_n} S_{m}$, the best estimator $\sh
_{\bar m}$ in a polynomial collection can be recovered by CV provided
$p/n$ converges to 1 as $n$ tends to $+\infty$.
The proof establishes this rate (i) cannot exceed $1/n$ to allow
distinguishing between parametric estimators (with convergence rate of
order $1/n$), and (ii) has to be faster than $(\log n)^{-1}$ to allow
dealing with the polynomial complexity of the model collection.
For instance, a finite collection would lead to replace the $(\log
n)^{-1}$ rate by a slower one determined by the control level of $\P
[ \mh= \bar m ]$.
In the regression setting, [\citet{Yang07}] already proved requiring
$p/n\to1$ enables to recover the best parametric estimator among
parametric ones (see Corollary~1), while this requirement is no longer
necessary when comparing parametric and nonparametric estimators.
Our result is consistent with Yang's one, although our setting is
somewhat different since we compare the best parametric estimator with
both parametric and nonparametric ones in the same time.

Conversely, when $s$ does not belong to $\bigcup_m S_{m}$, every parametric
model is biased according to \eqref{hyp.PNP.parametric} and $\sh
_{\bar
m}$ reaches a nonparametric rate, that is $nR_n(\bar m) \to+\infty$ as
$n$ tends to $+\infty$.
%
\begin{thmm}[(Model consistency with $s \notin\bigcup_m S_{m}$)] 
\label{thmm.model.consistency.not.belong}
Let $\bigcup_{m\in\M_n} S_{m}$ denote a collection of models
satisfying~\eqref{hyp.poly.collection} and~\eqref{hyp.param.nonparam},
$\bar m \in\mathcal{M}_n$ given by~\eqref
{expr.true.model.identification} be such
that~\eqref{hyp.best.model.identification} holds true, and
assume~\eqref
{hyp.sq.integrable.density},~\eqref
{hyp.regular.basis.dimension},~\eqref
{hyp.dimension.maximal} and~\eqref{hyp.expectation.lowerbound.dimension}.
For every $1\leq p\leq n-1$, let us also define $\mh= \mh(p) =\break
\argmin
_{m\in\M_n} \Rh_p(m)$.
Let us assume the target $s \notin\bigcup_{m\in\M_n} S_{m}$ and
%
%
$R_n(\bar m)\to0$ as $n$ tends to $+\infty$.
\begin{longlist}[1.]
\item[1.] If for large enough values of $n$ $D_{\bar m} \leq(\log n)^4$,
then every $1\leq p=p_n \leq n-1$ such that
%
\renewcommand{\theequation}{\arabic{equation}}
\setcounter{equation}{24}
\begin{equation}
\label{const.identification.nonparam.small} \log(n) \biggl( 1- \frac{p}{n} \biggr)\mathop{\longrightarrow}\limits_{n\to+
\infty} 0 \quad\mbox{and}\quad n\E\bigl[ \| s_{\bar m} - \sh_{\bar m}
\|^2 \bigr] = o( n-p )
\end{equation}
leads to
\[
\P[ \mh= \bar m ] \mathop{\longrightarrow}\limits_{n\to+\infty} 1.
\]

\item If for large enough values of $n$ $D_{\bar m} > (\log n)^4$, then
every $1\leq p=p_n \leq n-1$ such that
%
\begin{eqnarray}
\label{const.identification.nonparam.large} \frac{ (\log n)^5}{n\E[ \| s_{\bar m} - \sh_{\bar m} \|^2 ]
} &=& o\biggl( \frac{ p/n }{1- p/n} \biggr)
\quad\mbox{and}
\nonumber
\\[-8pt]
\\[-8pt]
\nonumber
\frac{p/n}{
1- p/n } &=& o\biggl( 1 \vee\frac{ \| s - s_{\bar m} \|^2 }{ \E [ \| s_{\bar
m} - \sh_{\bar m} \|^2 ] } \biggr)
\end{eqnarray}
leads to
\[
\P[ \mh= \bar m ] \mathop{\longrightarrow}\limits_{n\to+\infty} 1.
\]
\end{longlist}
\end{thmm}
The proof is similar to that of Theorem~\ref
{thmm.model.consistency.belong} and has been postponed to Section~C.1
[supplementary material \citet{Celi:2014:supp}].
The constraints on $p$ strongly depend on the rate of convergence of
$\sh_{\bar m}$ (nonparametric here).
When $S_{\bar m}$ is a small nonparametric model ($D_{\bar m}\leq(\log
n)^4$), \eqref{const.identification.nonparam.small} is very similar to
\eqref{conditions.p.belong} in the parametric setting. In particular,
$n\E[ \| s_{\bar m} - \sh_{\bar m} \|^2 ]\to+\infty$ as $n$
tends to $+\infty$ implies $n(1-p/n) \to+\infty$ as well.
For large nonparametric models ($D_{\bar m}> (\log n)^4$), the
constraints on $p$ are related to the ratio $\| s - s_{\bar m} \|^2/
\E[ \| s_{\bar m} - \sh_{\bar m} \|^2 ] $. For instance, when
estimating $s \in\mathcal{H}(L,\alpha)$ by regular histograms, this
ratio remains bounded while $n\E[ \| s_{\bar m} - \sh _{\bar m} \|^2
] $ grows polynomially in $n$. Then $p/n$ has to converge to $0$
as $n$ increases, but not too fast. In particular, Loo ($p=1$) is
suboptimal in that setting [see Figure~\ref
{fig.identification.notbelong} panel (b)].
Note that Theorem~\ref{thmm.model.consistency.not.belong} has the same
flavor as Corollary~1 in [\citet{Yang07}], except the density estimation
setting allows to relate $p$ to the features of the best estimator more closely.
%

%

\subsubsection{Simulation experiments}\label{subsubsec.Simuls.Identification}

\begin{figure}[t]

\includegraphics{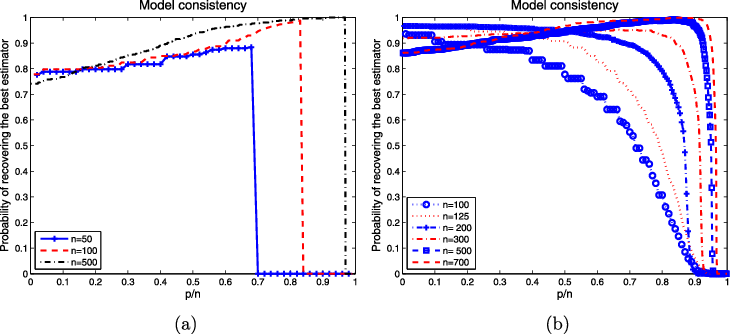}

\caption{$p/n \mapsto\P[ \mh= \bar m ]$ for density $s_1$ [panel
\textup{(a)}] and $s_2$ [panel \textup{(b)}]. $N=1
000$ samples have been drawn.} \label{fig.identification.belong}
\end{figure}

\begin{figure}

\includegraphics{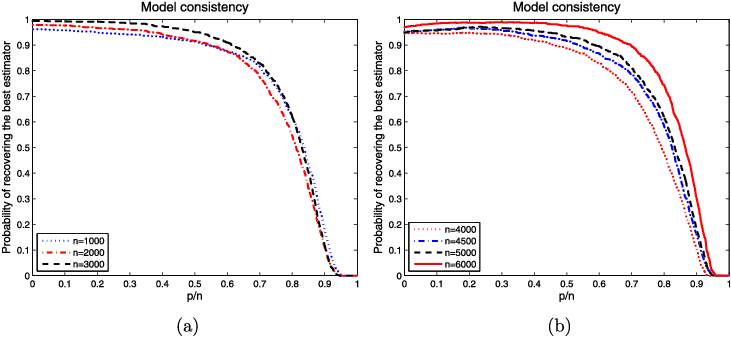}

\caption{$p/n \mapsto\P[ \mh = \bar m ]$ for density $s_3$ [panel
\textup{(a)}] and $s_4$ [panel \textup{(b)}].
$N=1 000$ samples have been drawn.} \label{fig.identification.notbelong}
\end{figure}

Simulation experiments have been performed in the settings of
Theorems~\ref{thmm.model.consistency.belong} and~\ref
{thmm.model.consistency.not.belong}, respectively, when $s$ belongs to
(resp., does not belong to) the model collection.
We used a polynomial model collection made of regular piecewise
constant functions described in Section~\ref
{subsubsec.projection.estimators} for which Assumptions~\eqref
{hyp.param.nonparam} and~\eqref{hyp.best.model.identification} are
fulfilled with $\mu=5.10^{-1}$.
In each setting, $N=1 000$ samples have been drawn.
Results are given in Figures~\ref{fig.identification.belong} and~\ref
{fig.identification.notbelong} where $\P[ \mh= \bar m ]$ is
displayed with respect to the ratio $p/n$.
Let us also mention that Lemma~B.3 in the supplement \citet
{Celi:2014:supp} clearly shows \eqref
{hyp.expectation.lowerbound.dimension} holds true for all densities
defined in the following as long as $D_m\geq2$ for every model in the
collection.

%

%

When \emph{$s$ belongs to the model collection} (Figure~\ref
{fig.identification.belong}), the following densities have been used:
\begin{enumerate}
\item$s_1(t) = \frac{6}{8} \1_{[0,1/2]}(t) + \frac{10}{8} \1_{[1/2,
1]}(t)$, $ t\in[0,1]$ [panel (a)],

\item$s_2(t) = \frac{135}{112} \1_{[0,1/3]}(t) + \frac{135}{56} \1
_{[1/3,1/2]}(t) + \frac{1}{4} \1_{[1/2, 5/7]}(t) + \frac{1}{2} \1_{[
5/7, 1]}(t)$, $ t\in[0,1]$ [panel (b)].
\end{enumerate}
As predicted by Theorem~\ref{thmm.model.consistency.belong}, CV reaches
model selection consistency for recovering the best parametric
estimator $\sh_{\bar m}$ on condition $p/n$ increases to $1$ as $n$
grows to $+\infty$.
Comparing (a) and (b), the convergence rate is slower in (b).
Unlike (a) where $\bar m$ remains almost unchanged as $n$ increases,
the best parametric estimator in (b) changes with $n$ as allowed
by~\eqref{expr.true.model.identification}.
Therefore, the slower convergence rate in~(b) results from the higher
dimension of the space of piecewise constant functions $s_2$ belongs to.

%

When \emph{$s$ does not belong to the model collection} (Figure~\ref
{fig.identification.notbelong}), densities with different smoothness
assumptions have been considered:
\begin{enumerate}
\item$s_3(x) = \beta(10,7; x)$, for every $t\in[0,1]$ [panel (a)],

\item$ s_4(x) = \frac{6}{5} x^{1/5}$, for every $t\in[0,1]$ [panel (b)].
\end{enumerate}
The converse situation arises since CV reaches model selection
consistency as long as $p/n$ decreases to $0$ as $n$ tends to $+\infty$.
Consistently with Theorem~\ref{thmm.model.consistency.not.belong}, this
rate strongly depends on the risk of the best estimator, that is on the
smoothness of the target.
While model selection consistency is illustrated by both panels (a)
and~(b), it is faster for the smoothest density $s_3$ than for $s_4$.
In Figure~\ref{fig.identification.notbelong} panel (b), the highest
probability is achieved for $p/n\approx0.18 $ with $n=6 000$.

%

\section{Discussion}
\label{sec.discussion}

From the present analysis of CV procedures in the density estimation
framework, we were able to prove the optimality of leave-one-out
cross-validation for risk estimation, which is consistent with earlier
results in the regression setting \citet{Burm89}.

However, when CV is used as model selection procedure, the optimal $p$
strongly depends on the structure of the model collection and on our
goal (estimation or identification).
%

\textit{Estimation}.
When the best model has dimension growing with $n$ [faster than $(\log
n)^{a}$ for some $a>0$] and the model collection has a polynomial
complexity \eqref{hyp.poly.collection}, Theorem~\ref
{thmm.poly.coll.regularity} proves any $p$ such that $p/n\to0$ leads to
an asymptotically optimal model selection procedure.
This is consistent with the asymptotic equivalence between Lpo (as long
as $p/n\to0$) and Mallows' $C_p$ previously settled in the regression
setting [\citet{Sha:1997}].

From a nonasymptotic point of view, Corollary~\ref
{cor.poly.coll.regularity} suggests choosing $p>1$ (for finite sample)
could balance the overfitting phenomenon arising from selecting a model
from a large collection.
This overfitting phenomenon is already well known with penalized
criteria such as Mallow's ones, inducing the need for heavier constants
in front of the penalty \citet{ArMa09}. Therefore, increasing $p$
amounts to penalize more strongly complex models (with large dimension).\looseness=1

\textit{Identification}.
As settled by \citet{Yang07} for regression, Section~\ref
{subsec.identification} highlights the optimal $p$ depends on the rate
of convergence of the estimator one tries to recover (and on the
structure of the model collection).

When the target estimator has a parametric rate with a polynomial
collection, Theorem~\ref{thmm.model.consistency.belong} proves $p/n\to
1$ leads to model selection consistency. This fact has been already
noticed by \citet{Shao93} in the regression setting who proved
leave-one-out is not model selection consistent.
Remembering the asymptotic equivalence between Lpo and BIC-like
criteria [\citet{Sha:1997}] established with the linear regression
model, this confirms the somewhat paradoxical requirement
[\citet{Yang06}] to devote most of available data ($p>n/2$) to the test set
when trying to recover a parametric estimator.

Drawing such a simple conclusion is harder when the best estimator has
a nonparametric rate as detailed by Theorem~\ref
{thmm.model.consistency.not.belong}.
If the best estimator has a rate close to parametric, then $p/n\to1$
provides model selection consistency.
Conversely, if the rate is slower (e.g., polynomial of order
$n^{-a}$, for some $a>0$), then requiring $p/n\to0$ enables to recover
the target estimator.
Relating that way the optimal $p$ to the rate of convergence of the
best estimator has been already done in the regression context by
[\citet{Yang07}, see Corollary~1].

Note that when the best estimator is nonparametric (e.g., with
polynomial rate), Theorem~\ref{thmm.poly.coll.regularity} and
Theorem~\ref{thmm.model.consistency.not.belong} imply $p/n \to0$ leads
to both efficiency and, respectively, model selection consistency.
However, there is no contradiction with the earlier paper by \citet
{Yang05} where it was proved no model selection criterion can share
both efficiency and model selection consistency in a parametric setting.
For instance, \citet{Li87} has established model selection consistency
for leave-one-out with nonparametric estimators in regression.\vadjust{\goodbreak}


\begin{appendix}\label{app}


\section{Estimation point of view}

\subsection{Main proofs}\label{app.proof.oracle}
\mbox{}
\begin{pf*}{Proof of Theorem~\ref{thmm.poly.coll.regularity}}
First, let\vspace*{1pt} us use Proposition~A.2 from \citet{Celi:2014:supp} applied
with $m,m'\in\M_n$ such that $ \widehat{R}_p(m') \leq\widehat{R}_p(m)$.
Then it comes
\begin{eqnarray*}
&& \frac{n}{n-p} \E\bigl[ Z_{m'}^2 \bigr] + \|
s-s_{m'} \|^2 - K(n,p) \bigl[ Z_{m'}^2
- \E\bigl[ Z_{m'}^2 \bigr] \bigr]
\\
&&\qquad \leq\frac{n}{n-p} \E\bigl[ Z_{m}^2 \bigr] + \|
s-s_{m} \|^2 - K(n,p) \bigl[ Z_{m}^2
- \E\bigl[ Z_{m}^2 \bigr] \bigr]
\\
&&\qquad\quad{} -2K(n,p) \nu_n( s_{m'} - s_m ) +
\frac{1}{n} \biggl( {K(n,p)+\frac
{n}{n-p}} \biggr) \nu_n(
\phi_{m'}-\phi_m ),
\end{eqnarray*}
where $K(n,p) = 1 + \frac{2}{n-1} + \frac{p}{n-p} \frac{1}{n-1}$.

Then, combining Propositions~B.4 and~B.5 from \citet{Celi:2014:supp} to
control the remainder terms, there exist a sequence $(\delta_n)_{\N}$
with $\delta_n\to0$ and $n\delta_n \to+\infty$ as $n\to+\infty$ and
an event $\Omega= \Omega_{\mathrm{rem},1} \cap\Omega_{\mathrm
{rem},2}$ of probability $1 - 4/n^2$ on which
\begin{eqnarray*}
& &\frac{n}{n-p} \E\bigl[ Z_{m'}^2 \bigr] + \|
s-s_{m'} \|^2 - K(n,p) \bigl[ Z_{m'}^2
- \E\bigl[ Z_{m'}^2 \bigr] \bigr]
\\
&&\qquad \leq\frac{n}{n-p} \E\bigl[ Z_{m}^2 \bigr] + \|
s-s_{m} \|^2 - K(n,p) \bigl[ Z_{m}^2
- \E\bigl[ Z_{m}^2 \bigr] \bigr]
\\
&&\qquad\quad{} + \delta_n K(n,p) \bigl( \| s-s_{m'} \|^2 +
\E\bigl[ Z_{m'}^2 \bigr] + \| s-s_{m} \|
^2 + \E\bigl[ Z_{m}^2 \bigr] \bigr)
\\
&&\qquad\quad{} + \delta_n \biggl( {K(n,p)+\frac{n}{n-p}} \biggr) \bigl[ \E
\bigl[ Z_{m'}^2 \bigr]+ \E\bigl[ Z_{m}^2
\bigr] \bigr].
\end{eqnarray*}
In the following, $\delta_n$ always denotes such a sequence even if the
precise expression of $\delta_n$ can differ from line to line.

Let us now use concentration results stated in Corollaries~B.1 and~B.2
from \citet{Celi:2014:supp} on the events $\Omega_{\mathrm{left}}$ and
$\Omega_{\mathrm{right}}$.
The important point in this proof is given by Lemmas~B.1 and~B
[\citet{Celi:2014:supp}], where it is proved that on
the event $\Omega=\Omega_{\mathrm{left}}\cap\Omega_{\mathrm
{right}}\cap\Omega_{\mathrm{rem},1}\cap\Omega_{\mathrm{rem},1}$,
$\min
\{ D_{m^*},D_{\mh(p)} \}\geq(\log n)^4$ for large enough values of $n$.
Therefore, one can apply Lemma~B.6 and Corollaries~B.1 and~B.2 from
\citet{Celi:2014:supp} with $L_m=0=r_n(m)$ to get
\begin{eqnarray*}
&& Z_{m'}^2 \biggl[ \biggl( \frac{n}{n-p} (1-
\delta_n) - 2\delta_n K(n,p) \biggr) ( 1-4
\delta_n ) - 4 K(n,p) \delta_n \biggr] \\
&&\quad{}+ \bigl[ 1 - \delta
_n K(n,p) \bigr] \| s-s_{m'} \|^2
\\
&&\qquad\leq Z_m^2 \biggl[ \biggl( \frac{n}{n-p}(1 +
\delta_n) + 2 \delta_n K(n,p) \biggr) ( 1+4
\delta_n ) + 4 K(n,p) \delta_n \biggr]\\
&&\qquad\quad{}+ \bigl[ 1 + \delta
_n K(n,p) \bigr] \| s-s_{m} \|^2.
\end{eqnarray*}

Choosing $m'=\mh$, it comes
\[
T_{V}^- Z_{\mh}^2 + T_{B}^- \|
s-s_{\mh} \|^2 \leq T_{V}^+ Z_{m}^2
+ T_{B}^+ \| s-s_{m} \|^2,
\]
where
\begin{eqnarray*}
T_{B}^- & = &1 - \delta_n K(n,p),\\
 T_{V}^- &=&
\frac
{n}{n-p}(1-\delta_n)[ 1 - 4 \delta_n ] - 2
K(n,p) \bigl[ 3\delta _n -4\delta_n^2 \bigr],
\\
T_{B}^+ & =& 1 + \delta_n K(n,p),\\
 T_{V}^+ &=&
\frac{n}{n-p}(1+\delta_n)[ 1 + 4 \delta_n ] + 2
K(n,p) \bigl[ 3\delta_n + 4\delta_n^2 \bigr].
\end{eqnarray*}

Finally on the event $\Omega$, the following oracle inequality holds
true for every $p\in\{ 1,n-1 \}$:
\[
\| s-\sh_{\mh(p)} \|^2 \leq C_{n}(p) \inf
_{m\in\mathcal{M}_n} \bigl\{ \| s-\sh _{m} \|^2 \bigr
\}\qquad \mbox{with } C_{n}(p) = \frac{ T_{B}^+ \vee
T_{V}^+ }{ T_{B}^- \wedge T_{V}^- }.
\]

Moreover, on the event $\Omega$, Lemmas~B.1 and~B.2 [\citet
{Celi:2014:supp}] show $\min\{ D_{m^*},D_{\mh(p)} \}\geq(\log n)^4$.
Then, it is enough to apply Propositions~B.1 and~B.2 from \citet
{Celi:2014:supp} to models satisfying this constraint, which leads to
the new event $\widetilde{\Omega}$ [where models with dimension smaller
than $(\log n)^4$ have been omitted] of probability at least $1-6/n^2$.
\end{pf*}

\begin{pf*}{Proof of Corollary~\ref{cor.poly.coll.regularity}}
Let us recall the expression of the leading constant
\[
C_{n}(p) = \frac{ T_{B}^+ \vee T_{V}^+ }{ T_{B}^- \wedge T_{V}^- },
\]
with
\begin{eqnarray*}
T_{B}^- & = &1 - \delta_n K(n,p), \\
 T_{V}^- &=&
\frac
{1}{1-p/n}(1-\delta_n)[ 1 - 4 \delta_n ] - 2
\delta_n K(n,p) [ 3 - 4\delta_n ],
\\
T_{B}^+ & = &1 + \delta_n K(n,p),\\
 T_{V}^+&=&
\frac{1}{1-p/n}(1+\delta_n)[ 1 + 4 \delta_n ] + 2\delta
_n K(n,p) [ 3 + 4\delta_n ],
\end{eqnarray*}
and $K(n,p) = 1 + \frac{2}{n-1} + \frac{p}{n-p} \frac{1}{n-1}$.

First, as long as $n$ is large enough, simple calculations when $p=1$
show $T_V^-(1)\leq T_B^-(1)$.
Noticing moreover that $T_V^+(p)\geq T_B^+(p)$ for every $p$, it comes
for $p$ close to 1
\[
C_{n}(p) = \frac{T_{V}^+ }{ T_{V}^- }  = \frac{ (1+\delta_n)[ 1 + 4
\delta_n ] + 2(1-p/n)\delta_n K(n,p) [ 3 + 4\delta_n ] }{
(1-\delta_n)[ 1 - 4 \delta_n ] - 2(1-p/n)\delta_n K(n,p)
[ 3 - 4\delta_n ] }.
\]
It is then easy to show that $p\mapsto C_n(p)$ is decreasing on $\{
1,\ldots,p^* \}$, where
$p^*$ denotes the value of $p$ such that $T_V^-(p)=T_B^-(p)$.
Hence,
\[
\frac{p_n^*}{n} = 1 - \frac{1 - 5\delta_n +4\delta_n^2 -( {2}/{(n-1)})
(3\delta_n - 4\delta_n^2) + {\delta_n}/{(n-1)} }{1 + 2(1+
{1}/{(n-1)}) (3\delta_n - 4\delta_n^2) - \delta_n(1+{1}/{(n-1)}) }.
\]
It results that for every $p\geq p^*$
\[
C_{n}(p) = \frac{T_{V}^+ }{ T_{B}^- },
\]
which is increasing with respect to $p$.

In the same way, it is easy to check that $p_n^*/ (10n \delta_n)
\mathop{\longrightarrow}\limits_{n\to+\infty}1 $, which enables us
to conclude the proof.
\end{pf*}

\begin{pf*}{Proof of Theorem~\ref{thmm.oracle.inequality.expectation}}
Introducing the event $\widetilde{\Omega}$ of Theorem~\ref
{thmm.poly.coll.regularity}, we get
\[
\E\bigl[ \| \bayes- \sh_{\mh(p)} \|^2 \bigr]  = \E\bigl[ \|
\bayes- \sh_{\mh
(p)} \|^2 \1_{\widetilde{\Omega}} \bigr] + \E\bigl[
\| \bayes- \sh_{\mh(p)} \|^2 \1_{\widetilde{\Omega}^c} \bigr].
\]
Then Theorem~\ref{thmm.poly.coll.regularity} applied to the first
expectation in the right-hand side leads to
\[
\E\bigl[ \| \bayes- \sh_{\mh(p)} \|^2 \bigr]  \leq
C_n(p) \E \Bigl[ \inf_{m\in\mathcal{M}_n} \| \bayes-
\sh_{m} \|^2 \Bigr] + \E\bigl[ \| \bayes -
\sh_{\mh(p)} \|^2 \1_{\widetilde{\Omega}^c} \bigr].
\]
Applying \eqref{hyp.bias}, one gets
\[
\E\bigl[ \| \bayes- s_{\mh(p)} \|^2 \1_{\widetilde{\Omega }^c} \bigr]
\leq\E\biggl[ \frac{ c_u }{ D_{\mh(p)}^u } \1_{\widetilde{\Omega}^c} \biggr] \leq
c_u \P\bigl( \widetilde{\Omega}^c \bigr) \leq
\frac{6 c_u}{n^2},
\]
while \eqref{hyp.regular.basis.dimension} and \eqref
{hyp.dimension.maximal} provide
\begin{eqnarray*}
&&\E\bigl[ \| s_{\mh(p)} - \sh_{\mh(p)} \|^2
\1_{\widetilde {\Omega }^c} \bigr]\\
&&\qquad = \E\biggl[ \sum_{\lambda\in\Lambda(\mh(p))} (
P_n\varphi_{\lambda
} - P\varphi_{\lambda} )^2
\1_{\widetilde{\Omega}^c} \biggr]
\\
& &\qquad\leq 2 \E\biggl[ \sum_{\lambda\in\Lambda(\mh(p))} ( P_n
\varphi _{\lambda} )^2 \1_{\widetilde{\Omega}^c} \biggr] + 2 \E\biggl[
\sum_{\lambda
\in \Lambda (\mh(p))} ( P\varphi_{\lambda} )^2
\1_{\widetilde
{\Omega}^c} \biggr]
\\
&&\qquad \leq 2 \E\Biggl[ \sum_{\lambda\in\Lambda(\mh(p))} \frac {1}{n^2}
\sum_{i, j=1}^n \varphi_{\lambda}(X_i)
\varphi_{\lambda}(X_j) \1 _{\widetilde{\Omega}^c} \Biggr] + 2 \| s
\|^2 \E[ D_{\mh(p)} \1_{\widetilde{\Omega}^c} ]
\\
&&\qquad \leq 2 \bigl( \Phi+ \| s \|^2\bigr) \frac{n}{(\log n)^2}\P\bigl(
\widetilde {\Omega}^c \bigr) \leq\bigl( \Phi+ \| s \|^2
\bigr) \frac{12}{n (\log n)^2}.
\end{eqnarray*}
\upqed\end{pf*}

\section{Identification point of view}\label{subsec.appendix.identification}

\subsection{Proof of Theorem~\texorpdfstring{\protect\ref{thmm.model.consistency.belong}}{3.3}}\label
{subsec.main.proofs.identification}
%
The general purpose is to prove there exist an event $\Omega_n$ with
$\P
( \Omega_n )\to1$ as $n$ tends to $+\infty$ and a positive integer
$N$ such that on~$\Omega_n$, for every $n\geq N$, every $m\neq\bar m$
satisfies
%
\renewcommand{\theequation}{\arabic{equation}}
\setcounter{equation}{26}
\begin{equation}
\label{ineq.ref.model.consistency.first} \Rh_p(\bar m) - \Rh_p(m) \leq-
u_n(m) \bigl(1 + o(1)\bigr),
\end{equation}
where $u_n(m)>0$ denotes a real number for every $n$ and $m$, and
$o(1)$ does not depend on $m$.
In particular, this implies
%
\begin{equation}
\label{conv.identification.modelconsistency} \P( \mh= \bar m ) = \P\bigl( \forall m\neq\bar m, \Rh
_p(\bar m) - \Rh _p(m) <0 \bigr) \mathop{\longrightarrow}\limits_{n\to+
\infty} 1,
\end{equation}
which would complete the proof.

Let us consider the event $\Omega_{\mathrm{left}}\cap\Omega
_{\mathrm
{right}}$ in Proposition~C.2 from \citet{Celi:2014:supp} with $\beta
_1=\beta_2=1/n^2$, and the events $\Omega_{\mathrm{rem},1}$
[Proposition~B.4 from \citet{Celi:2014:supp}] and $\Omega_{\mathrm
{rem},3}$ (Proposition~C.1 [\citet{Celi:2014:supp}]).
Then with $\Omega_n = \Omega_{\mathrm{left}} \cap\Omega_{\mathrm
{right}} \cap\Omega_{\mathrm{rem},1} \cap\Omega_{\mathrm{rem},3}$
and $\P[ \Omega_n^c ] \leq8/n^2$, showing \eqref
{conv.identification.modelconsistency} amounts to prove
\[
\P\bigl( \Omega_n\cap\bigl\{ \forall m\neq\bar m,
\Rh_p(\bar m) - \Rh_p(m) <0 \bigr\} \bigr)
\mathop{\longrightarrow}\limits_{n\to+\infty} 1.
\]

Let us now focus on the event $\Omega_n$.
The two main steps correspond to distinguishing between \emph
{parametric} and \emph{nonparametric} models $S_m$ [see \eqref
{hyp.param.nonparam}].
For every $m$, let us define $B(m) = \| s - s_m \|^2$ and $V(m) = \E
[ \| s_m - \sh_m \|^2 ]$, where $s_m = \argmin_{t\in S_m}
\| s - t \|^2$.
From line to line, the value of $\delta_n$ may change, but it always
denotes a sequence decreasing to $0$ and such that $n \delta_n \to
+\infty$ as $n$ grows.

\textit{If $\sh_m$ has a parametric rate}.

\begin{itemize}
\item If $s\in S_{m}$:

Let us first notice $s_m=s=s_{\bar m}$, which implies $R_n(m)=V(m)$ and
$R_n(\bar m)=V(\bar m)$.
Then
Proposition~A.2 from \citet{Celi:2014:supp}, and Propositions~B.4
and~C.2 from \citet{Celi:2014:supp} lead to
\begin{eqnarray*}
&&\biggl\vert\bigl[ \Rh_p(\bar m) - \Rh_p(m) \bigr] -
\frac{n}{n-p} \bigl[ R_n(\bar m) - R_n(m) \bigr]
\biggr\vert \\
&&\qquad\leq \biggl( 36 L_n^2 + 3 \delta_n +
\delta_n \frac{n}{n-p} \biggr) \bigl[ R_n(m) +
R_n(\bar m) \bigr]
\\
&&\qquad=  o\biggl( \frac{n}{n-p} \biggr) \bigl[ R_n(m) +
R_n(\bar m) \bigr],
\end{eqnarray*}
by requiring $L_n^2=o( (1-p/n)^{-1} )$, which provides \eqref
{ineq.ref.model.consistency.first} by use of \eqref
{hyp.best.model.identification}.
Note that in the previous inequality, $ r_n(\bar m)$ and $r_n(m)$ [from
Proposition~C.2 from \citet{Celi:2014:supp}] have been omitted since
they are negligible with respect to the other terms.

\item If $s\notin S_{m}$:

Similarly, Proposition~A.2, Proposition~B.4, and Propositions~C.1
and~C.2 from \citet{Celi:2014:supp} lead to
\begin{eqnarray*}
&&\biggl \vert\bigl[ \Rh_p(\bar m) - \Rh_p(m) \bigr] - \bigl[
B(\bar m) - B(m) \bigr] - \frac
{n}{n-p} \bigl[ V(\bar m) - V(m) \bigr] \biggr\vert
\\
&&\qquad\leq 6\delta_n \bigl[ B(\bar m) + B(m) \bigr]\\
&&\qquad\quad{} + \biggl(
36 L_n^2 + 3\delta _n \frac
{\Phi}{\xi} +
\delta_n \biggl( 3 + \frac{n}{n-p} \biggr) \biggr) \bigl[ V(m) +
V(\bar m) \bigr]
\\
&&\qquad\quad{} + 3\delta_n \| s \| \sqrt{\frac{\Phi}{\xi}} \sqrt{V(\bar m) +
V(m) } + r_n(\bar m).
\end{eqnarray*}
With $s\in\bigcup_{m'}S_{m'}$ and $s\notin S_m$, it comes $B(\bar m)=0$
and $B(m)\geq\tau>0$ by \eqref{hyp.param.nonparam}.
Since both $\widehat{s}_m$ and $\sh_{\bar m}$ have parametric rates, requiring
$n(1-p/n) \to+\infty$ as $n$ grows implies \eqref
{ineq.ref.model.consistency.first}, that is,
\[
\bigl[ \Rh_p(\bar m) - \Rh_p(m) \bigr] \leq- B(m)
\bigl( 1 + o(1) \bigr).
\]
\end{itemize}

\textit{If $\sh_m$ has a nonparametric rate}.

\begin{itemize}
\item If $s\in S_m$:

Proposition~A.2 from \citet{Celi:2014:supp},
$s_m = s = s_{\bar m}$, and Propositions~C.1 and~C.2 from \citet
{Celi:2014:supp} combined with $L_n^2= o( n/(n-p) )$ provide
\begin{eqnarray*}
 &&\biggl\vert\bigl[ \Rh_p(\bar m) - \Rh_p(m) \bigr] -
\frac{n}{n-p} \bigl[ V(\bar m) - V(m) \bigr]\biggr \vert \\
&&\qquad\leq o\biggl(
\frac{n}{n-p} \biggr) V(m) + o\biggl( \frac {n}{n-p} \biggr) V(\bar m),
\end{eqnarray*}
where $o( n/(n-p) ) $ does not depend on $m$.
Since $S_m$ is nonparametric, \eqref{hyp.param.nonparam} gives
\[
\frac{ V(\bar m) }{ V(m) } \leq\frac{\pi}{ (\log n)^2 } \Bigl( n/(\log n)^2 \inf
_m V(m) \Bigr)^{-1} \mathop{\longrightarrow}\limits_{n\to+\infty} 0,
\]
which implies \eqref{ineq.ref.model.consistency.first} with
\[
(n-p)\bigl[ \Rh_p(\bar m) - \Rh_p(m) \bigr] \leq- n
V(m) \bigl( 1 + o(1) \bigr).
\]

\item If $s\notin S_m$:

Both $L_n^2= o( n/(n-p) )$ and the same argument as above lead to
\begin{eqnarray*}
&&\biggl\vert\bigl[ \Rh_p(\bar m) - \Rh_p(m) \bigr] + B(m) +
\frac{n}{n-p} V(m) \bigl[ 1 + o(1) \bigr]\biggr \vert \\
&&\qquad\leq V(m) o\biggl(
\frac{n}{n-p} \biggr) + B(m) o( 1 ).
\end{eqnarray*}
Then \eqref{ineq.ref.model.consistency.first} holds true with
\[
\bigl[ \Rh_p(\bar m) - \Rh_p(m) \bigr] \leq- B(m)
\bigl(1+o(1)\bigr) - \frac{n}{n-p} V(m) \bigl(1+o(1)\bigr).
\]
\end{itemize}
Then there exists an integer $N$ such that for $n\geq N$, on the event
$\Omega_n$, \eqref{ineq.ref.model.consistency.first} holds true for
every $m\in\M_n$, which completes the proof.
\end{appendix}

\section*{Acknowledgements}
We thank Sylvain Arlot and St\'ephane Robin for helpful discussions,
and also the Associate Editor as well as two anonymous referees for
their comments and suggestions to improve the earlier version of the paper.


\begin{supplement}
\stitle{Supplement to ``Optimal cross-validation in density estimation
with the $L^2$-loss'': Technical proofs and details}
\slink[doi]{10.1214/14-AOS1240SUPP} 
\sdatatype{.pdf}
\sfilename{aos1240\_supp.pdf}
\sdescription{Owing to space constraints, we have moved technical
proofs to a supplementary document [\citet{Celi:2014:supp}].}
\end{supplement}


%

\printaddresses
\end{document}